\documentclass[preprint,3p,sort&compress,final,times]{elsarticle}
\usepackage{amsmath}
\usepackage{amssymb}
\usepackage{mathtools}
\usepackage{placeins}
\usepackage[usenames,dvipsnames]{color}
\usepackage{xcolor}
\usepackage{cases}
\usepackage{array}

\usepackage[usenames,dvipsnames]{color}

\newcommand{\thethanks}{This work was supported by computational resources provided by the Australian Government through NCI under the National Computational Merit Allocation Scheme (NCMAS). We thank the APR.Intern program and the Australian Mathematical Sciences Institute for their support and funding.}

\journal{}
\begin{document}
\begin{frontmatter}

\title{A comparison of variational upwinding schemes for geophysical fluids, and their application
	to potential enstrophy conserving discretisations}
\author[BOM]{David Lee\corref{cor}}
\ead{davelee2804@gmail.com}
\author[ANU]{Alberto F. Mart\'{\i}n}
\author[BOM]{Christopher Bladwell}
\author[MON,CIMNE]{Santiago Badia}
\address[BOM]{Bureau of Meteorology, Melbourne, Australia.}
\address[ANU]{School of Computing, Australian National University, Acton ACT 2601, Australia.}
\address[MON]{School of Mathematics, Monash University, Melbourne 3800, Australia.}
\address[CIMNE]{Centre Internacional de M\`{e}todes Num\`{e}rics en Enginyeria, Esteve Terrades 5, E-08860 Castelldefels, Spain.}
\cortext[cor]{Corresponding author. Tel. +61 452 262 804.}

\begin{abstract}
Methods for stabilising the turbulent cascade of potential enstrophy are analysed and compared for a compatible finite element discretisation of the rotating shallow water equations. 
These different approaches to upwinding the potential vorticity include the well-known anticipated potential vorticity method (APVM), streamwise upwind Petrov-Galerkin (SUPG) method, and a more
recent method where the trial functions are evaluated downstream within the reference element. In all cases the upwinding scheme
conserves both potential vorticity and energy, since the antisymmetric structure of the equations is preserved. The APVM leads to a
symmetric definite correction to the potential enstrophy that is dissipative and inconsistent.
The SUPG scheme introduces a consistent correction to the APVM scheme that acts as a backscatter term ensuring a richer depiction of
turbulent dynamics.
The downwinded trial function formulation results in the advection of downwind corrections, which analysis and numerical experiments
show to be quantitatively similar to the SUPG scheme for a turbulent shear flow.
The main difference between the SUPG and downwinded trial function schemes is in the energy conservation and residual errors. If just two nonlinear iterations
are applied then the energy conservation errors are improved for the downwinded trial function formulation, reflecting a smaller
residual error than for the SUPG scheme.

We also present new temporal formulations by which potential enstrophy is exactly integrated across each time level. Results using these
formulations are observed to be stable in the absence of any dissipation, despite the aliasing of grid scale turbulence.
Using such a formulation and the APVM with a coefficient $\mathcal{O}(100)$ times smaller
that its regular value leads to turbulent spectra that are greatly improved at the grid scale over the SUPG and downwinded 
trial function formulations with unstable potential enstrophy errors.
\end{abstract}

\end{frontmatter}

\section{Introduction}

Two dimensional turbulence in geophysical fluids involves a nonlinear cascade of energy to large scales and
of potential enstrophy to small scales. Consequently such flows may be stabilised via the removal of potential
enstrophy at small scales, via the upwinding of the potential vorticity, without breaking energy conservation.
This may be achieved by replacing the potential vorticity in the vector invariant form of the momentum equation
with a value sampled at some upstream location, as is the case for the well-known anticipated potential vorticity
method (APVM) \cite{SB85}, or alternatively via a variational approach by which the test functions are weighted
towards an upstream location, such as the streamwise upwind Petrov-Galerkin (SUPG) \cite{BH82} or variational
multiscale method (VMS) \cite{Hughes98}. 
Both the APVM \cite{MC14} and SUPG \cite{NC17} methods have previously been used
in finite element models of geophysical systems for the energetically consistent transport of material quantities
(vorticity \cite{NC17}, potential vorticity \cite{Natale16,Bauer18} and potential temperature {\cite{WCB21}}). A
similar approach has also been applied in a finite volume scheme for the shallow water equations derived from a 
variational principle \cite{Brecht21}, and VMS approaches have been applied to spectral element discretisations
of the compressible Euler equations for atmospheric dynamics \cite{MG15}.

In the present study the APVM and SUPG methods are compared to a more recent
approach by which the potential vorticity trial functions are evaluated at downstream locations within the reference
element \cite{Lee21}
in terms of their consistency, conservation properties and residual errors. This reference element stabilisation
has been used previously for the stabilisation of both the potential enstrophy for the shallow water equations,
in which case it is applied by downwinding the trial functions \cite{Lee21}, and to the potential
temperature in the 3D compressible Euler equations, in which case it is applied as the upwinding of the test
functions \cite{LP21}. In the context of strong form differential operators, for which the differential operator
is applied to the trial function (such as the divergence of a mass flux in $H(div,\Omega)$ for the flux form
advection of some quantity in $L^2(\Omega)$), this is achieved via the upwinding of the test function (provided
that the test function is $C^0$ continuous across element boundaries in the direction of the flow). For a weak form differential operator
however (such as material advection as the adjoint of the aforementioned flux form operator) this reference element
upwinding is conversely applied via the \emph{downwinding} of the trial functions \cite{Lee21}.

The principal contributions of this paper are the comparison of the recently developed downwinded trial function
method of potential enstrophy stabilisation to the well established APVM and SUPG methods, both in terms of
analysis and numerical experiments. We also introduce a new method of exactly integrating the potential 
enstrophy across each time level. For an implicit time integrator this new temporal representation of the potential 
enstrophy allows for the simulation of well developed turbulence on the sphere without any form of dissipation whatsoever, since
the turbulent cascade of potential entrophy to the grid scale is well controlled.

The remainder of this article proceeds as follows: In Section 2, the mixed finite element discretisation of the
rotating shallow water equations and its conservation properties are introduced. The APVM, SUPG and downwinded
trial function formulations of the potential enstrophy are analysed in the context of this discretisation
in Section 3. Section 3 also introduces new time discretisations for the potential vorticity at first and second 
order, which allow for the exact integration of the potential enstrophy across the time level, such that 
simulations may be run to a mature turbulent state without any form of dissipation.
Results comparing
the three upwinding schemes for the solution of a shear flow instability test case \cite{Galewsky04} are presented in
Section 4, and conclusions are discussed in Section 5.

\section{Mixed finite element discretisation the rotating shallow water equations}

The rotating shallow water equations for the evolution of the velocity, $\boldsymbol{u}$ and depth, $h$,
within a two dimensional domain $\Omega\subset\mathbb{R}^2$
may be expressed as
\begin{subequations}\label{eq::sw}
\begin{align}
	\frac{\partial\boldsymbol{u}}{\partial t} + q\boldsymbol{F}^{\perp} + \nabla P &= 0,\\
\frac{\partial h}{\partial t} + \nabla\cdot\boldsymbol{F} &= 0,
\end{align}
\end{subequations}
where
\begin{subequations}\label{eq::var_derivs}
\begin{align}
	\boldsymbol{F} &= h\boldsymbol{u},\\
	P &= \frac{1}{2}\boldsymbol{u}\cdot\boldsymbol{u}+gh,
\end{align}
\end{subequations}
are the mass flux, $\boldsymbol{F}$ and the Bernoulli potential, $P$ respectively,
with $g$ being the acceleration due to gravity, and 
\begin{equation}\label{eq::pv}
q = \frac{\nabla\times\boldsymbol{u}+f}{h},
\end{equation}
is the potential vorticity, while $f$ is the Coriolis term.
The operator $\perp$ denotes rotation as 
$\boldsymbol{F}^{\perp} = [-\boldsymbol{F}\cdot\hat{\boldsymbol{y}},\boldsymbol{F}\cdot\hat{\boldsymbol{x}}]$
for the unit vectors $[\hat{\boldsymbol{x}}, \hat{\boldsymbol{y}}]$ which specify a coordinate system for $\Omega$.
We further note that both $q$ and $f$ are scalar quantities within $\Omega$, however for the extension to a three 
dimensional domain these may be regarded as the radial components of a vector quantity.

In order to spatially discretise the above system we introduce finite dimensional subspaces
$\mathbb{V}_0\subset H^1(\Omega)$, $\mathbb{V}_1\subset H(div,\Omega)$, $\mathbb{V}_2\subset L^2(\Omega)$, which
may be compatibly mapped to one another via the discrete De-Rham complex of the form \cite{MC14,LPG18,Eldred18}
\begin{equation}\label{H_rot_H_div_L_2}
\mathbb{R}^2\longrightarrow\mathbb{V}_0 \stackrel{\nabla^{\perp}}{\longrightarrow}\mathbb{V}_1
\stackrel{\nabla\cdot}{\longrightarrow}\mathbb{V}_2\longrightarrow 0.
\end{equation}
The prognostic equations \eqref{eq::sw} are then solved for the discrete approximations $q_h\in\mathbb{V}_0$,
$\boldsymbol{u}_h,\boldsymbol{F}_h\in\mathbb{V}_1$, $h_h,P_h\in\mathbb{V}_2$ by multiplying by the test functions
$\boldsymbol{v}_h\in\mathbb{V}_1$, $\phi_h\in\mathbb{V}_2$ respectively. 
We emphasise the difference here between \emph{trial} functions, basis functions used
to represent the solution fields, $q_h$, $\boldsymbol{u}_h$, $\boldsymbol{F}_h$, $h_h$, $P_h$, 
and \emph{test} functions, $\boldsymbol{v}_h$, $\phi_h$,
basis functions against which the discrete form of the differential equations
are weighted in order to derive a variational form of the system.
Integrating over the
domain $\Omega$ using integrating by parts and assuming periodic boundary conditions, 
the spatially discrete form of the equations are given as
\begin{subequations}\label{eq::sw_disc}
\begin{align}
	\int\boldsymbol{v}_h\cdot\frac{\partial\boldsymbol{u}_h}{\partial t}\mathrm{d}\Omega +
	\int\boldsymbol{v}_hq_h\boldsymbol{F}_h^{\perp}\mathrm{d}\Omega -
	\int\nabla\cdot\boldsymbol{v}_h P_h\mathrm{d}\Omega &= 0\qquad\forall\boldsymbol{v}_h\in\mathbb{V}_1,\label{eq::mom_eqn_disc}\\
	\int\phi_h\frac{\partial h_h}{\partial t}\mathrm{d}\Omega +
	\int\phi_h\nabla\cdot\boldsymbol{F}_h\mathrm{d}\Omega &= 0\qquad\forall\phi_h\in\mathbb{V}_2.\label{eq::cont_eqn_disc}
\end{align}
\end{subequations}

The discrete form of the energy is given as
\begin{equation}
	\mathcal{H}_h(\boldsymbol{u}_h,h_h) = \int\frac{1}{2}h_h\boldsymbol{u}_h\cdot\boldsymbol{u}_h + \frac{1}{2}gh_h^2\mathrm{d}\Omega.
\end{equation}
The variational derivatives of the energy with respect to the solution variables, $\boldsymbol{u}_h$, $h_h$
may be computed discretely as the perturbation of the solution variables in some arbitrary test direction 
$\boldsymbol{v}_h$, $\phi_h$ by some small value $\epsilon$ in the limit that $\epsilon\rightarrow 0$
\cite{Celledoni12}.
The corresponding mass flux and Bernoulli potential are therefore given in analogy to \eqref{eq::var_derivs} as
\begin{subequations}\label{eq::var_derivs_disc}
\begin{align}
	\frac{\mathrm{d}}{\mathrm{d}\epsilon}\mathcal{H}_h(\boldsymbol{u}_h+\epsilon\boldsymbol{v}_h,h_h)\Big|_{\epsilon=0} &= 
	\int\boldsymbol{v}_h\cdot\frac{\delta\mathcal{H}_h}{\delta\boldsymbol{u}_h}\mathrm{d}\Omega = 
	\int\boldsymbol{v}_h\cdot\boldsymbol{F}_h\mathrm{d}\Omega =
	\int\boldsymbol{v}_h\cdot\boldsymbol{u}_hh_h\mathrm{d}\Omega\qquad\forall\boldsymbol{v}_h\in\mathbb{V}_1,\\
	\frac{\mathrm{d}}{\mathrm{d}\epsilon}\mathcal{H}_h(\boldsymbol{u}_h,h_h+\epsilon\phi_h)\Big|_{\epsilon=0} &= 
	\int\phi_h\frac{\delta\mathcal{H}_h}{\delta h_h}\mathrm{d}\Omega = 
	\int\phi_h P_h\mathrm{d}\Omega =
	\int\phi_h\Bigg(\frac{1}{2}\boldsymbol{u}_h\cdot\boldsymbol{u}_h + gh_h\Bigg)\mathrm{d}\Omega\qquad\forall\phi_h\in\mathbb{V}_2.
\end{align}
\end{subequations}
Moreover, the discrete form of the potential vorticity, $q_h\in\mathbb{V}_0$ is determined in analogy to \eqref{eq::pv} as
\begin{equation}\label{eq::pv_disc}
	\int\psi_hh_hq_h\mathrm{d}\Omega = \int-\nabla^{\perp}\psi_h\cdot\boldsymbol{u}_h + \psi_hf_h\mathrm{d}\Omega
	\qquad\forall\psi_h\in\mathbb{V}_0,
\end{equation}
with $f_h\in\mathbb{V}_0$ being the discrete representation of the Coriolis term.

\subsection{Conservation properties}

In order to describe the semi-implicit solution of the rotating shallow water equations, we first introduce the 
dual spaces $\mathbb{V}_0'$, $\mathbb{V}_1'$ and $\mathbb{V}_2'$
and the $L^2(\Omega)$ inner product as 
$(a,b) = \int ab\mathrm{d}\Omega$. The duality pairing is given as $\langle a,\cdot\rangle$, which maps coefficients of
$a$ to coefficients in the corresponding dual space. The linear mass matrices are then given as
$\boldsymbol{\mathsf{M}}_1:\mathbb{V}_1\rightarrow\mathbb{V}_1'$ for which 
$\boldsymbol{\mathsf{M}}_1\boldsymbol{v}_h := \langle\boldsymbol{v}_h,\cdot\rangle$, $\forall\boldsymbol{v}_h\in\mathbb{V}_1$ and 
$\boldsymbol{\mathsf{M}}_2:\mathbb{V}_2\rightarrow\mathbb{V}_2'$ for which 
$\boldsymbol{\mathsf{M}}_2\phi_h := \langle\phi_h,\cdot\rangle$, $\forall\phi_h\in\mathbb{V}_2$. 
We also introduce the divergence operator,
$\boldsymbol{\mathsf{D}}:\mathbb{V}_1\rightarrow\mathbb{V}_2'$ as 
$\boldsymbol{\mathsf{D}}\boldsymbol{v}_h:=\langle\nabla\cdot\boldsymbol{v}_h,\cdot\rangle$ and the rotational operator
$\boldsymbol{\mathsf{C}}:\mathbb{V}_0\times\mathbb{V}_1\rightarrow\mathbb{V}_1'$ as 
$\boldsymbol{\mathsf{C}}(\psi_h,\boldsymbol{v}_h):=\langle\psi_h\boldsymbol{v}_h^{\perp},\cdot\rangle$.

The energy conserving structure of the spatially and temporally discrete formulation of the rotating
shallow water equations \eqref{eq::sw_disc} and \eqref{eq::var_derivs_disc} is exposed by
expressing these in skew-symmetric form as:
\begin{equation}\label{eq::sw_ss}
	\begin{bmatrix}
		\boldsymbol{\mathsf{M}}_1(\boldsymbol{u}_h^{n+1}-\boldsymbol{u}_h^n) \\
		\boldsymbol{\mathsf{M}}_2(h_h^{n+1}-h_h^n)
	\end{bmatrix} + \Delta t
	\begin{bmatrix}
		\boldsymbol{\mathsf{C}}(\cdot,\overline{q}_h) &
		-\boldsymbol{\mathsf{D}}^{\top} \\
		\boldsymbol{\mathsf{D}} & \boldsymbol{\mathsf{0}}
	\end{bmatrix}
	\begin{bmatrix}
		\overline{\boldsymbol{F}}_h \\ \overline{P}_h
	\end{bmatrix} = 0
\end{equation}
where $\overline{\boldsymbol{F}}_h$, $\overline{P}_h$ are the exact time integrals of the 
mass flux and Bernoulli potential respectively.
Since we compute solutions only at the beginning and end of the time level using our
implicit time integrator, and not at intermediate times, we can only assume solutions that 
are piecewise linear in time. These piecewise linear representations of the variational 
derivatives can be integrated exactly between time levels $n$ and $n+1$ using Simpson's rule,
resulting in implicit expressions for the exact second order temporal integrals of the variational 
derivatives \eqref{eq::var_derivs_disc} as \cite{Hairer11,Bauer18,Lee21b}
\begin{subequations}\label{eq::temp_var_derivs}
	\begin{align}
		\int\boldsymbol{v}_h\cdot\overline{\boldsymbol{F}}_h\mathrm{d}\Omega &=
		\frac{1}{6}\int\boldsymbol{v}_h\cdot\boldsymbol{u}_h^{n+1}(2h_h^{n+1}+h_h^n) +
			       \boldsymbol{v}_h\cdot\boldsymbol{u}_h^n(h_h^{n+1}+2h_h^n)
			       \mathrm{d}\Omega\qquad\forall\boldsymbol{v}_h\in\mathbb{V}_1,\\
		\int\phi_h\overline{ P}_h\mathrm{d}\Omega &= \label{eq::dH_dh}
		\frac{1}{6}\int\phi_h(\boldsymbol{u}_h^{n+1}\cdot\boldsymbol{u}_h^{n+1} +
				      \boldsymbol{u}_h^{n+1}\cdot\boldsymbol{u}_h^n +
				      \boldsymbol{u}_h^n\cdot\boldsymbol{u}_h^n)\mathrm{d}\Omega +
		\frac{g}{2}\int\phi_h(h_h^{n+1}+h_h^n)
			       \mathrm{d}\Omega\qquad\forall\phi_h\in\mathbb{V}_2.
	\end{align}
\end{subequations}
Note that the variational derivatives are continuous both \emph{within} and \emph{between} time
levels, since $\boldsymbol{u}_h^{n+1}$, $h_h^{n+1}$ at the end of a given time level $n$ are the same as
the values at the beginning of the next time level, $n+1$.

Energy conservation is then
achieved for the discrete form by setting the test functions as the variational derivatives:
$\boldsymbol{v}_h = \overline{\boldsymbol{F}}_h$, $\phi_h=\overline{ P}_h$ and
adding the two equations together. Cancellation of the right hand sides results in the expression
\begin{equation}\label{eq::dHdt}
	\Bigg(\frac{\delta\mathcal{H}_h}{\delta\boldsymbol{u}_h},\frac{(\boldsymbol{u}_h^{n+1} - \boldsymbol{u}_h^n)}{\Delta t}\Bigg) +
	\Bigg(\frac{\delta\mathcal{H}_h}{\delta h_h},\frac{(h_h^{n+1}-h_h^n)}{\Delta t}\Bigg) =
	\frac{\mathrm{d}\mathcal{H}_h}{\mathrm{d}t} = 0.
\end{equation}
The above relation holds due to the skew-symmetry of the spatial discretisation, and the exact temporal
integral of the variational derivatives, such that the temporal chain rule above is preserved discretely 
\cite{LP21}.

In addition to the energy, the spatially discrete shallow water equations may also conserve potential enstrophy,
\begin{equation}
	\mathcal{Z}_h = 
	\int\frac{h_hq_h^2}{2}\mathrm{d}\Omega.
\end{equation}
This is achieved by recalling the compatible mappings supported
by the De-Rham complex in \eqref{H_rot_H_div_L_2} and setting
$\boldsymbol{v}_h=\nabla^{\perp}\psi_h$ in \eqref{eq::mom_eqn_disc} (for $\psi_h\in\mathbb{V}_0$ and 
$\boldsymbol{v}_h\in\mathbb{V}_1$) \cite{MC14}. Since the De-Rham complex further
ensures that $\nabla\cdot\nabla^{\perp}:= 0$ is preserved discretely, and recalling \eqref{eq::pv_disc},
(using the fact that $f_h$ does not vary with time) and \eqref{eq::cont_eqn_disc} this implies a 
potential vorticity conservation equation (for periodic boundary conditions) of the form \cite{MC14}
\begin{equation}\label{eq::pv_adv}
	\int\psi_h\frac{\partial h_hq_h}{\partial t}\mathrm{d}\Omega + \int\psi_h\nabla\cdot(q_h\boldsymbol{F}_h)\mathrm{d}\Omega = 0
	\qquad\forall\psi_h\in\mathbb{V}_0.
\end{equation}
Since the Coriolis term is independent of time, the vorticity, $\omega_h = h_hq_h - f_h$ is conserved by
the above expression for a choice of $\psi_h=1$.

Potential enstrophy conservation may be derived by instead setting $\psi_h=q_h$ in \eqref{eq::pv_adv},
and recalling the discrete continuity equation \eqref{eq::cont_eqn_disc} \cite{MC14,LPG18}. However this may
be more generally derived by regarding the potential enstrophy as an additional
invariant of the Hamiltonian system, or \emph{Casimir}, for which conservation is ensured due to the fact
that its variational derivatives are part of a null space for the skew-symmetric operator in \eqref{eq::sw_ss}.
To see this we first compute the variational derivative with respect to the velocity as
\begin{equation}\label{eq::var_deriv_Z_u}
	\int\boldsymbol{v}_h\cdot\frac{\delta\mathcal{Z}_h}{\delta\boldsymbol{u}_h}\mathrm{d}\Omega =
	\Bigg(\frac{\delta q_h}{\delta\boldsymbol{u}_h}\Bigg)^{\top}\frac{\delta\mathcal{Z}_h}{\delta q_h}.
\end{equation}
In order to describe these variational derivatives we introduce the additional linear operators 
$\boldsymbol{\mathsf{R}}:\mathbb{V}_0\rightarrow\mathbb{V}_1'$ as
$\boldsymbol{\mathsf{R}}\psi_h:=\langle\nabla^{\perp}\psi_h,\cdot\rangle$ and 
$\boldsymbol{\mathsf{H}}_0(\phi_h,\psi_h):\mathbb{V}_0\otimes\mathbb{V}_2\rightarrow\mathbb{V}_0'$ as
$\boldsymbol{\mathsf{H}}_0(\phi_h,\psi_h) :=\langle \phi_h\psi_h,\cdot\rangle$.
Following \eqref{eq::pv_disc}, the variational derivative of the potential vorticity $q_h$ with respect to the velocity
$\boldsymbol{u}_h$ is given as the matrix 
\begin{equation}\label{eq::var_deriv_q_u}
	\int \psi_h h_h\frac{\delta q_h}{\delta\boldsymbol{u}_h}\mathrm{d}\Omega = -\int\nabla^{\perp}\psi_h \cdot \boldsymbol{u}_h\mathrm{d}\Omega
	\Longrightarrow\frac{\delta q_h}{\delta\boldsymbol{u}_h} = -\boldsymbol{\mathsf{H}}^{-1}_0(h_h,\boldsymbol{\mathsf{R}}^{\top} \boldsymbol{u}_h)
\end{equation}
and the variational derivative of the functional $\mathcal{Z}_h$ with respect to $q_h$ returns the vector
\begin{equation}\label{eq::var_deriv_Z_q}
	\int\psi_h\frac{\delta\mathcal{Z}_h}{\delta q_h}\mathrm{d}\Omega = \int\psi_hh_hq_h\mathrm{d}\Omega
	\Longrightarrow  \frac{\delta\mathcal{Z}_h}{\delta q_h} = \boldsymbol{\mathsf{M}}_0^{-1} \boldsymbol{\mathsf{H}}_0(h_h,q_h).
\end{equation}
Combining \eqref{eq::var_deriv_q_u}, \eqref{eq::var_deriv_Z_q} into \eqref{eq::var_deriv_Z_u}, the 
variational derivative of $\mathcal{Z}_h$ with respect to $\boldsymbol{u}_h$ in $\mathbb{V}_1$ is given as
\begin{equation}
	\frac{\delta\mathcal{Z}_h}{\delta\boldsymbol{u}_h} =
	-\boldsymbol{\mathsf{M}}_1^{-1}\boldsymbol{\mathsf{R}}\boldsymbol{\mathsf{H}}^{-1}_0(h_h,\boldsymbol{\mathsf{H}}_0(h_h,q_h)) = 
	-\boldsymbol{\mathsf{M}}_1^{-1}\boldsymbol{\mathsf{R}}q_h.
\end{equation}
Introducing the operator
$\boldsymbol{\mathsf{Q}}:\mathbb{V}_0\times\mathbb{V}_0\rightarrow\mathbb{V}_2'$ as
$\boldsymbol{\mathsf{Q}}(\psi_h,\chi_h):=\langle\psi_h\chi_h,\cdot\rangle$,
the variational derivative of the potential enstrophy with respect to the fluid
depth in $\mathbb{V}_2$ is similarly given as
\begin{equation}
	\frac{\delta\mathcal{Z}_h}{\delta h_h} = \frac{1}{2}\boldsymbol{\mathsf{M}}^{-1}_2\boldsymbol{\mathsf{Q}}(q_h,q_h).
\end{equation}
Multiplying the skew-symmetric operator in \eqref{eq::sw_ss} by the vector
\begin{equation}
	\left[\Bigg(\frac{\delta\mathcal{Z}_h}{\delta\boldsymbol{u}_h}\Bigg)^{\top},\quad
	\Bigg(\frac{\delta\mathcal{Z}_h}{\delta h_h}\Bigg)^{\top}\right] = 
	\left[
		-q_h\boldsymbol{\mathsf{R}}^{\top}\boldsymbol{\mathsf{M}}_1^{-1},\quad
		\frac{1}{2}(q_h,q_h)\boldsymbol{\mathsf{Q}}^{\top}\boldsymbol{\mathsf{M}}^{-1}_2
	\right]
\end{equation}
then gives
\begin{subequations}
\begin{align}
	-q_h\boldsymbol{\mathsf{R}}^{\top}\boldsymbol{\mathsf{M}}_1^{-1}&
	\Bigg(\boldsymbol{\mathsf{C}}(\cdot,q_h) - \boldsymbol{\mathsf{D}}^{\top}\Bigg) +
	\frac{1}{2}(q_h,q_h)\boldsymbol{\mathsf{Q}}^{\top}\boldsymbol{\mathsf{M}}^{-1}_2\boldsymbol{\mathsf{D}} \\
	&= -\nabla^{\perp}q_h\boldsymbol{\mathsf{C}}(\cdot,q_h) +
	\Bigg(\boldsymbol{\mathsf{D}}\boldsymbol{\mathsf{M}}_1^{-1}\boldsymbol{\mathsf{R}}q_h\Bigg)^{\top} + 
	\frac{q_h^2}{2}\boldsymbol{\mathsf{D}}\\
	&= \int q_h\boldsymbol{v}_h\cdot\nabla q_h\mathrm{d}\Omega + \int\phi_h\nabla\cdot\nabla^{\perp}q_h\mathrm{d}\Omega +
	      \int\frac{q_h^2}{2}\nabla\cdot\boldsymbol{v}_h\mathrm{d}\Omega\\
	&= \int\nabla\Bigg(\frac{q_h^2}{2}\Bigg)\cdot\boldsymbol{v}_h\mathrm{d}\Omega +
	      \int\frac{q_h^2}{2}\nabla\cdot\boldsymbol{v}_h\mathrm{d}\Omega =
	      \int\nabla\cdot\Bigg(\frac{q_h^2}{2}\boldsymbol{v}_h\Bigg)\mathrm{d}\Omega = 0.
\end{align}
\end{subequations}
Where in the last line we have recalled the discrete property of the compatible bases that $\nabla\cdot\nabla^{\perp}:=0$,
as well as the assumption of periodic boundary conditions.
Note that the last line above is contingent on the discrete preservation of the chain rule as
\begin{equation}
	\int q_h\nabla(q_h)\cdot\boldsymbol{v}_h\mathrm{d}\Omega :=
	\int\frac{1}{2}\nabla(q_h^2)\cdot\boldsymbol{v}_h\mathrm{d}\Omega,
\end{equation}
which only holds for exact integration \cite{LPG18}, and so is not strictly preserved for a trigonometric
Jacobian on the surface of the sphere involving transcendental functions. Now that the spatial terms have been
shown to cancel for the discrete form of the potential enstrophy evolution equation, we are left with the
temporal chain rule expansion in analogy to \eqref{eq::dHdt} as
\begin{equation}\label{eq::pe_adv}
	\Bigg(\frac{\delta\mathcal{Z}_h}{\delta\boldsymbol{u}_h}\frac{(\boldsymbol{u}_h^{n+1}-\boldsymbol{u}_h^n)}{\Delta t}\Bigg) +
	\Bigg(\frac{\delta\mathcal{Z}_h}{\delta h_h}\frac{(h_h^{n+1}-h_h^n)}{\Delta t}\Bigg) =
	\frac{\mathrm{d}\mathcal{Z}_h}{\mathrm{d}t} = \frac{\mathrm{d}}{\mathrm{d}t}\int\frac{h_hq_h^2}{2}\mathrm{d}\Omega = 0.
\end{equation}

\subsection{Exact time integration of potential enstrophy}

In this section new temporal formulations of the potential vorticity are introduced by which the potential enstrophy
is integrated exactly across the time levels. Doing so leads to enhanced stability, such that we are able to
run our simulations of two dimensional turbulence at low Mach number without stabilisation.

Since we do not explicitly solve the advection equation for the potential vorticity \eqref{eq::pv_adv}, there 
is no enforcement of continuity of $q_h$ between time levels. Rather this is diagnosed from
\eqref{eq::pv_disc}.
Integrating both sides of \eqref{eq::pv_disc} with respect to $q_h$ gives
\begin{equation}
	\int h_hq_h^2\mathrm{d}\Omega =
	-2\int\nabla^{\perp}q_h\cdot\boldsymbol{u}_h\mathrm{d}\Omega + 2\int q_hf_h\mathrm{d}\Omega.
\end{equation}
The second order integral of the above expression across the time level from $t^n$ to $t^{n+1}$ for 
piecewise linear in time $q_h$, $\boldsymbol{u}_h$ and $h_h$
yields a cubic temporal polynomial
on the left hand side, and a quadratic on the right. Both of these may be integrated exactly using Simpson's
rule, as
\begin{subequations}
	\begin{align}
		2\mathcal{Z}_h :=&
	\frac{1}{6}\int\frac{3}{2}h_h^n(q_h^n)^2 + \frac{1}{2}h_h^n(q_h^{n+1})^2 + (h_h^n+h_h^{n+1})q_h^nq_h^{n+1} +
		\frac{1}{2}h_h^{n+1}(q_h^{n})^2 + \frac{3}{2}h_h^{n+1}(q_h^{n+1})^2\mathrm{d}\Omega \\ =&
	-\frac{1}{3}\int \nabla^{\perp}q_h^n\cdot(2\boldsymbol{u}_h^n + \boldsymbol{u}_h^{n+1}) +
	\nabla^{\perp}q_h^{n+1}\cdot(\boldsymbol{u}_h^n + 2\boldsymbol{u}_h^{n+1})\mathrm{d}\Omega +
	\int (q_h^n+q_h^{n+1})f_h\mathrm{d}\Omega.
	\end{align}
\end{subequations}
Differentiating with respect to $q_h^n$ and $q_h^{n+1}$ leads to the coupled system for the solution of the
potential vorticity at both time levels
\begin{equation}\label{eq::q_linear_in_time}
\begin{bmatrix}
	\frac{1}{6}\boldsymbol{\mathsf{H}}_0(\cdot,3h_h^n + h_h^{n+1}) & \frac{1}{6}\boldsymbol{\mathsf{H}}_0(\cdot,h_h^n + h_h^{n+1}) \\
	\frac{1}{6}\boldsymbol{\mathsf{H}}_0(\cdot,h_h^n + h_h^{n+1}) & \frac{1}{6}\boldsymbol{\mathsf{H}}_0(\cdot,h_h^n + 3h_h^{n+1})
\end{bmatrix}
\begin{bmatrix}q_h^n \\ q_h^{n+1}\end{bmatrix} = 
	\begin{bmatrix}
		-\frac{1}{3}\boldsymbol{\mathsf{R}}^{\top}(2\boldsymbol{u}_h^n + \boldsymbol{u}_h^{n+1}) +
		\boldsymbol{\mathsf{M}}_0f_h \\
		-\frac{1}{3}\boldsymbol{\mathsf{R}}^{\top}(\boldsymbol{u}_h^n + 2\boldsymbol{u}_h^{n+1}) +
		\boldsymbol{\mathsf{M}}_0f_h
	\end{bmatrix}
\end{equation}
where as before $\boldsymbol{\mathsf{M}}_0:\mathbb{V}_0\rightarrow\mathbb{V}_0'$ and is given as
$\boldsymbol{\mathsf{M}}_0\psi_h:=\langle\psi_h,\cdot\rangle$.
The above system is linear in $q_h^n$, $q_h^{n+1}$, and may be solved at each nonlinear iteration of
\eqref{eq::sw_ss}.
We then compute the mean potential vorticity across the time level simply as $\overline{q}_h=(q_h^n+q_h^{n+1})/2.$

Alternatively, we may represent $\overline{q}_h$ as being piecewise constant in time, for which the 
potential enstrophy is exactly integrated over the time level as
\begin{equation}
2\mathcal{Z}_h :=
	\frac{1}{2}\int (h_h^n + h_h^{n+1})\overline{q}_h^2\mathrm{d}\Omega =
	-\int \nabla^{\perp}\overline{q}_h\cdot(\boldsymbol{u}_h^n + \boldsymbol{u}_h^{n+1})\mathrm{d}\Omega +
	2\int\overline{q}_hf_h\mathrm{d}\Omega.
\end{equation}
In this case the potential vorticity is 
more simply diagnosed from the piecewise linear in time $\boldsymbol{u}_h$, $h_u$ by differentiating with 
respect to $\overline{q}_h$ as
\begin{equation}\label{eq::q_const_in_time}
	\boldsymbol{\mathsf{H}}_0(\overline{q}_h,h_h^n + h_h^{n+1}) = 
	-\boldsymbol{\mathsf{R}}^{\top}(\boldsymbol{u}_h^n + \boldsymbol{u}_h^{n+1}) +
	2\boldsymbol{\mathsf{M}}_0f_h.
\end{equation}
We stress that while the above formulations yield exact temporal integration of potential enstrophy across the
time level, they do not ensure exact conservation \emph{between} time levels, such that unlike the energy, 
potential enstrophy will still not be exactly conserved in time.

\subsection{Time integration}

The resulting nonlinear problem \eqref{eq::sw_ss}, \eqref{eq::var_derivs_disc} and \eqref{eq::q_linear_in_time} or 
\eqref{eq::q_const_in_time} is then solved using a Newton method with a constant in
time approximate Jacobian \cite{Bauer18,Wimmer20}, which at each nonlinear iteration $k$ is solved for the
updates $\delta\boldsymbol{u}_h^k = \boldsymbol{u}_h^{k+1} - \boldsymbol{u}_h^k$,
$\delta h_h^k = h_h^{k+1} - h_h^k$ as
\begin{equation}\label{eq::newton_system}
	\begin{bmatrix}
		\boldsymbol{\mathsf{M}}_1 + 0.5\Delta t\boldsymbol{\mathsf{C}}(\cdot,f) &
		-0.5\Delta tg\boldsymbol{\mathsf{D}}^{\top} \\
		0.5\Delta tH\boldsymbol{\mathsf{D}} & \boldsymbol{\mathsf{M}}_2
	\end{bmatrix}
	\begin{bmatrix}
		\delta\boldsymbol{u}_h^k \\ \delta h_h^k
	\end{bmatrix} = -
	\begin{bmatrix}
		R_{\boldsymbol{u}}^k \\ R_{h}^k
	\end{bmatrix}
\end{equation}
where $H$ is the mean fluid depth and $k$ is the estimate of the solution at time level $n+1$ 
at the current Newton iteration.

Note that this constant Jacobian omits nonlinear terms due to variations in the potential vorticity, mass flux
and Bernoulli potential. Since the vorticity field varies only slowly with time compared to the gravity waves,
and we consider only weakly hyperbolic flow regimes for which variations in $h_h$ are small with respect to $H$,
this is sufficient to achieve numerical convergence. The residuals at iteration $k$ are given 
via \eqref{eq::sw_ss} as:
\begin{subequations}\label{eq::residuals}
	\begin{align}
		R_{\boldsymbol{u}}^k &= \boldsymbol{\mathsf{M}}_1(\boldsymbol{u}_h^k - \boldsymbol{u}_h^n) +
		\Delta t\boldsymbol{\mathsf{C}}(\overline{\boldsymbol{F}}_h,\overline{q}_h) -
		\Delta t\boldsymbol{\mathsf{D}}^{\top}\overline{P}_h
		\label{eq::mom_resid}\\
		R_{h}^k &= \boldsymbol{\mathsf{M}}_2(h_h^k - h_h^n) +
		\Delta t\boldsymbol{\mathsf{D}}\overline{\boldsymbol{F}}_h.
	\end{align}
\end{subequations}
Note that at each Newton iteration $k$, prior to the evaluation of the solution increments
$\delta\boldsymbol{u}_h^k$, $\delta h_h^k$, we must first evaluate the variational derivatives, 
$\overline{\boldsymbol{F}}_h$, $\overline{P}_h$, as well as the potential vorticity $\overline{q}_h$
(setting $k=n+1$ as a proxy at each iteration). The full solution procedure at each nonlinear iteration
$k$ therefore follows as:
\begin{enumerate}
	\item Compute the exact integral of the variational derivatives $\overline{\boldsymbol{F}}_h$, $\overline{P}_h$ 
		across the time level $n$ and the approximation of time level $n+1$ at iteration $k$ via \eqref{eq::temp_var_derivs}.
	\item Compute the second order integral of the potential vorticity $q_h$ as either:
		\begin{enumerate}
			\item A time centered average $q_h^{n+1/2} = (q_h^n + q_h^k)/2$, where $q_h^n$ and $q_h^k$ are
				evaluated instantaneously via \eqref{eq::pv_disc}.
			\item Exact piecewise constant integration of $\overline{q}_h$ over the time level via \eqref{eq::q_const_in_time}.
			\item Exact piecewise linear integration of $\overline{q}_h$ over the time level via \eqref{eq::q_linear_in_time}.
		\end{enumerate}
	\item Evaluate the residuals $R_{\boldsymbol{u}}^k$, $R_h^k$ via \eqref{eq::residuals}
	\item Solve for the solution increments $\delta\boldsymbol{u}_h^k$, $\delta h_h^k$ via \eqref{eq::newton_system}
	\item Update the solutions as $\boldsymbol{u}_h^{k+1} = \boldsymbol{u}_h^{k} + \delta\boldsymbol{u}_h^k$,
		$h_h^{k+1} = h_h^{k} + \delta h_h^k$, and increment $k\rightarrow k+1$.
	\item Terminate the Newton iteration if $k > k_{max}$, for maximum iteration number $k_{max}$, or if 
		$|\delta\boldsymbol{u}_h^k|/|\boldsymbol{u}_h^k| < \epsilon_{tol}$ and 
		$|\delta h_h^k|/|h_h^k| < \epsilon_{tol}$, where $\epsilon_{tol}$ is some specified solver tolerance.
\end{enumerate}

\section{Potential vorticity upwinding}

\subsection{The anticipated potential vorticity method (APVM)}

Energy conservation is preserved for any choice of $q_h$ in \eqref{eq::mom_eqn_disc}, since the term
$\int\boldsymbol{v}_hq_h\boldsymbol{F}_h^{\perp}\mathrm{d}\Omega = 0$ for $\boldsymbol{v}_h=\boldsymbol{F}_h$
is in any case skew-symmetric, as shown in \eqref{eq::sw_ss}. Consequently $q_h$ may be upwinded as a
means of damping small scale oscillations without breaking energy conservation. Perhaps the most
straight-forward form of upwinding is the anticipated potential vorticity method \cite{SB85,MC14},
for which the potential vorticity is replaced in \eqref{eq::mom_eqn_disc} by its upstream value
\begin{equation}\label{eq::up_apvm}
	q_h-\tau\boldsymbol{u}_h\cdot\nabla q_h,
\end{equation}
where $\tau$ is an upwinding timescale parameter. In the simplest case, $\tau$ is typically a constant
of the same order as the time step, $\mathcal{O}(\Delta t$),
however it may also include additional terms that scale with the advective term or other physics \cite{AT04}.
Substitution of \eqref{eq::up_apvm} into \eqref{eq::mom_eqn_disc}, and then following the same procedure 
used to derive the flux form advection equation for $q_h$ in \eqref{eq::pv_adv} \cite{MC14},
results in a modified form of the potential vorticity advection equation \eqref{eq::pv_adv} as
\begin{equation}\label{eq::q_adv_apvm}
\int\psi_h\frac{\partial h_hq_h}{\partial t} +
\psi_h\nabla\cdot\Bigg(\boldsymbol{F}_h\Bigg(q_h - \tau\boldsymbol{u}_h\cdot\nabla q_h\Bigg)\Bigg)\mathrm{d}\Omega = 0
\qquad\forall\psi\in\mathbb{V}_0.
\end{equation}
Note that this expression is not consistent with respect to \eqref{eq::pv_adv}, such that it will converge
to a different solution for $q_h$ with increased resolution. Again setting $\psi_h=q_h$, the analogous form
of \eqref{eq::pe_adv} is given as
\begin{equation}\label{eq::pe_adv_apvm}
\int\frac{\partial}{\partial t}\Bigg(\frac{h_hq_h^2}{2}\Bigg)\mathrm{d}\Omega +
\int\tau h_h(\boldsymbol{u}_h\cdot\nabla q_h)^2\mathrm{d}\Omega = 0.
\end{equation}
From this we see that the APVM upwinding results in the addition of a symmetric positive definite dissipation
term to the potential enstrophy evolution that acts to always remove potential enstrophy from the system for
$\tau > 0$.

\subsection{The streamwise upwind Petrov-Galerkin method (SUPG)}

As a consistent alternative to the APVM \eqref{eq::up_apvm}, we may instead 
replace the potential vorticity by a corrected upwinded formulation as
\cite{Natale16,Bauer18}
\begin{equation}
	q_h-\tau\Bigg(\frac{\partial q_h}{\partial t} + \boldsymbol{u}_h\cdot\nabla q_h\Bigg),
\end{equation}
for which
\begin{equation}
	\int\psi_h\frac{\partial h_hq_h}{\partial t} + \psi_h\Bigg(\nabla\cdot\boldsymbol{F}_h\Bigg
	(q_h - \tau\Bigg(\frac{\partial q_h}{\partial t} + \boldsymbol{u}_h\cdot\nabla q_h\Bigg)\Bigg)\Bigg)\mathrm{d}\Omega = 0
	\qquad\forall\psi_h\in\mathbb{V}_0.
\end{equation}
Rearranging this expression, and again recalling the continuity equation \eqref{eq::cont_eqn_disc} gives
\begin{equation}
	\int h_h\Bigg(\psi_h + \tau\boldsymbol{u}_h\cdot\nabla \psi_h\Bigg)
	\Bigg(\frac{\partial q_h}{\partial t} + \boldsymbol{u}_h\cdot\nabla q_h\Bigg)\mathrm{d}\Omega = 0
	\qquad\forall\psi_h\in\mathbb{V}_0.
\end{equation}
This consistent modification of the trial function constitutes the SUPG formulation \cite{BH82,NC17}. The corresponding
potential enstrophy evolution equation is given for $\psi_h=q_h$ as
\begin{equation}\label{eq::pe_adv_supg}
	\int\frac{\partial}{\partial t}\Bigg(\frac{h_hq_h^2}{2}\Bigg)\mathrm{d}\Omega +
	\int\tau h_h(\boldsymbol{u}_h\cdot\nabla q_h)^2\mathrm{d}\Omega +
	\int\tau h_h
	\Bigg(\frac{\partial q_h}{\partial t}\Bigg)
	(\boldsymbol{u}_h\cdot\nabla q_h)
	\mathrm{d}\Omega = 0.
\end{equation}
The second term in \eqref{eq::pe_adv_supg} is identical to that in \eqref{eq::pe_adv_apvm}, ensuring that 
potential enstrophy is dissipated via the same mechanism as for APVM. Since $q_h$ is materially advected, 
we have that $\frac{q_h}{dt}\approx-\boldsymbol{u}_h\cdot\nabla q_h$, such that the third term acts as a 
backscatter term to inject potential enstrophy back into the system, ensuring a richer representation of 
grid scale turbulence than for the APVM. However since this term is not sign definite it will not fully 
counter balance the potential enstrophy dissipation of the second term, ensuring the net dissipation of 
the overall scheme.

Note that we have made no assumption as to the form of $\tau$, other than that it is $>0$. This may be set as a constant value, or
as a spatially and temporally varying field which scales with both the time derivative and the advection term.

\subsection{Downwinded trial functions}

As an alternative to these traditional approaches, one may effectively upwind the potential vorticity by
\emph{downwinding} the reference element trial functions in the case of weak form differential operators
\eqref{eq::pv_disc}, or \emph{upwinding} the test functions in the case of strong form differential operators
\cite{Lee21,LP21}.
In either case, this may be achieved by moving the location of the quadrature point, $\xi$
at which the reference element basis function, $\hat{l}_i(\xi)$ is evaluated, provided that
$\hat{l}_i(\xi); i=1,...,p$ span a finite element space of polynomial degree $p$ that is $C^0$
continuous in the direction of
the flow (which is the case in both dimensions for the space $\mathbb{V}_0\subset H^1(\Omega)$).

For the \emph{downwinded} reference element trial functions, these are evaluated as
$\hat{l}_h^d(\xi;\hat{u}(\xi)) = \hat{l}_h(\xi - \tau\hat{u}(\xi)/|J|)$,
where $\tau$ is again some time scale parameter, $\hat{u}$ is the velocity in reference element coordinates,
and $|J|$ is the Jacobian determinant of the geometrical map, $\Phi$ \cite{Rognes06}, evaluated at $\xi$.
Conversely, the \emph{upwinded} test functions
may be evaluated as $\hat{l}_h^u(\xi;\hat{u}(\xi)) = \hat{l}_h(\xi + \tau\hat{u}(\xi)/|J|)$. For a trial
(resp. test) function of degree $p$ may be expanded as
\begin{subequations}
\begin{align}
	\hat{l}_i^d(\xi;\hat{u}(\xi)) &= \hat{l}_i\Big(\xi - \frac{\tau\hat{u}(\xi)}{|J|}\Big) = \hat{l_i}(\xi) -
	\frac{\tau\hat{u}(\xi)}{|J|}\frac{d\hat{l}_i(\xi)}{d\xi} +
	\frac{(\tau\hat{u}(\xi))^2}{2|J|^2}\frac{d^2\hat{l}_i(\xi)}{d\xi^2} -
	\frac{(\tau\hat{u}(\xi))^3}{6|J|^3}\frac{d^3\hat{l}_i(\xi)}{d\xi^3} + ...
	\frac{(-\tau\hat{u}(\xi))^p}{p!|J|^p}\frac{d^p\hat{l}_i(\xi)}{d\xi^p}, \label{eq::basis_down}\\
	\hat{l}_i^u(\xi;\hat{u}(\xi)) &= \hat{l}_i\Big(\xi + \frac{\tau\hat{u}(\xi)}{|J|}\Big) = \hat{l_i}(\xi) +
	\frac{\tau\hat{u}(\xi)}{|J|}\frac{d\hat{l}_i(\xi)}{d\xi} +
	\frac{(\tau\hat{u}(\xi))^2}{2|J|^2}\frac{d^2\hat{l}_i(\xi)}{d\xi^2} +
	\frac{(\tau\hat{u}(\xi))^3}{6|J|^3}\frac{d^3\hat{l}_i(\xi)}{d\xi^3} + ...
	\frac{(+\tau\hat{u}(\xi))^p}{p!|J|^p}\frac{d^p\hat{l}_i(\xi)}{d\xi^p}.
\end{align}
\end{subequations}
In two dimensions the downwinded trial/upwinded test functions may be expressed in reference space as
\begin{subequations}
\begin{align}\label{eq::q_down}
	\hat{\psi}_h^d(\boldsymbol{\xi};\hat{\boldsymbol{u}}) &=
	\hat{l}_i^d(\xi;\hat{u}(\boldsymbol{\xi}))\otimes \hat{l}_j^d(\eta;\hat{v}(\boldsymbol{\xi})),\\
	\hat{\psi}_h^u(\boldsymbol{\xi};\hat{\boldsymbol{u}}) &=
	\hat{l}_i^u(\xi;\hat{u}(\boldsymbol{\xi}))\otimes \hat{l}_j^u(\eta;\hat{v}(\boldsymbol{\xi})),
\end{align}
\end{subequations}
where $\boldsymbol{\xi} = (\xi,\eta)$ are the two dimensional local element coordinates,
$\hat{\boldsymbol{u}} = (\hat{u},\hat{v})$ are the two dimensional
velocity components in the reference element and $h=j(p+1)+i$ is the index of the two dimensional tensor product
basis function.
These downwinded (resp. upwinded) trial (resp. test) functions may then be pushed forward into physical space as
\begin{subequations}
\begin{align}
	\psi_h^d(\boldsymbol{x};\hat{\boldsymbol{u}}(\boldsymbol{\xi})) &= 
	\left[\hat{\psi}_h\Bigg(\boldsymbol{\xi} - \frac{\tau\hat{\boldsymbol{u}}(\boldsymbol{\xi})}{|J|}\Bigg)\right]\circ\Phi^{-1}(\boldsymbol{x}),\\
	\psi_h^u(\boldsymbol{x};\hat{\boldsymbol{u}}(\boldsymbol{\xi})) &= 
	\left[\hat{\psi}_h\Bigg(\boldsymbol{\xi} + \frac{\tau\hat{\boldsymbol{u}}(\boldsymbol{\xi})}{|J|}\Bigg)\right]\circ\Phi^{-1}(\boldsymbol{x}).
\end{align}
\end{subequations}
where $\boldsymbol{x}$ are the global coordinates. The downwinded form of the potential enstrophy is then given in
physical coordinates by $q_h^d(\boldsymbol{x})$.

Since it is the trial functions that are downwinded, the above expansion enters both the flux term and the temporal
derivative, giving the analogue of \eqref{eq::pv_adv} as
\begin{equation}
	\int\psi_h\frac{\partial h_hq_h^d}{\partial t} +
	\psi_h\nabla\cdot(\boldsymbol{F}_hq_h^d)
	\mathrm{d}\Omega = 0\qquad\forall\psi_h\in\mathbb{V}_0,
\end{equation}
Since the downwinded correction to the potential vorticity in \eqref{eq::q_down} can be expressed as a Taylor series
correction to the static trial functions, we label this as $q_h' = q_h^d - q_h$. Rearranging and recalling the continuity
equation then gives
\begin{equation}\label{eq::up_trial}
	\int\psi_h h_h\Bigg(\frac{\partial q_h}{\partial t} + \boldsymbol{u}_h\cdot\nabla q_h\Bigg)\mathrm{d}\Omega
	+ \int\psi_h h_h
	\Bigg(\frac{\partial}{\partial t} + \boldsymbol{u}_h\cdot\nabla\Bigg)
	q_h'\mathrm{d}\Omega = 0\qquad\forall\psi_h\in\mathbb{V}_0.
\end{equation}
While the additional trial function corrections presented above are not consistent, if instead the test
functions are \emph{upwinded} using the same approach, as has been done previously for the stabilisation
of potential temperature in the 3D compressible Euler equations \cite{LP21}, then the resulting stabilisation
is indeed consistent.
The resulting potential enstrophy evolution equation for the downwinded trial functions is given as
\begin{equation}\label{eq::pe_adv_up_q}
	\int\frac{\partial}{\partial t}\Bigg(\frac{h_hq_h^2}{2}\Bigg)\mathrm{d}\Omega + \int h_hq_h
	\Bigg(\frac{\partial q_h'}{\partial t} + \boldsymbol{u}_h\cdot\nabla q_h'\Bigg)
	\mathrm{d}\Omega = 0.
\end{equation}
From \eqref{eq::basis_down} we have to a first approximation that $q'_h\approx-\tau\boldsymbol{u}_h\cdot\nabla q_h$.
If we further assume periodic boundary conditions in time as well as space for the grid scale dynamics driven
by the potential enstrophy cascade, then we may apply integration by parts to the material transport operator 
$\mathcal{L}=\frac{\partial}{\partial t} + \boldsymbol{u}_h\cdot\nabla$ in the second term of \eqref{eq::pe_adv_up_q},
resulting in an expression equivalent to \eqref{eq::pe_adv_supg}. We therefore anticipate that for the small
scales on which the damping term operates that the downwindind trial function formulation should yield similar
results as for the SUPG scheme, and any consistency errors will be small.

As stated above, the downwinded trial functions are applied both when the potential vorticity is diagnosed,
resulting in a downwind correction to the time derivative in the advection equation, and to the rotational
term in the momentum equation, resulting in a downwinded correction to the flux term in the advection
equation. The fully discrete diagnostic \eqref{eq::pv_disc} and residual \eqref{eq::mom_resid} expressions
are therefore expressed for the downwinded trial functions as
\begin{equation}
	\boldsymbol{\mathsf{H}}_0(q_h^d,h_h) = -\boldsymbol{\mathsf{R}}^{\top}\boldsymbol{u}_h + \boldsymbol{\mathsf{M}}_0f_h,
\end{equation}
\begin{equation}
	\boldsymbol{\mathsf{M}}_1R_{\boldsymbol{u}}^k = \boldsymbol{\mathsf{M}}_1(\boldsymbol{u}_h^k - \boldsymbol{u}_h^n) +
	\Delta t\boldsymbol{\mathsf{C}}(\overline{\boldsymbol{F}}_h,\overline{q}_h^d) -
	\Delta t\boldsymbol{\mathsf{D}}^{\top}\overline{P}_h.
\end{equation}

\section{Results}

\subsection{Nonlinearly converged solutions}

In order to quantify the differences between the APVM, SUPG and downwinded trial function formulations,
these are compared in the context of a shear flow instability of a barotropic jet, as triggered by a
small initial perturbation in the otherwise geostrophically balanced depth field \cite{Galewsky04},
on a cubed sphere using polynomials of degree $p=3$ for the representation of the potential vorticity
\cite{LP18}. The test case was run for 20 days using $6\times 32\times 32$ elements and 8 point
Gauss-Lobatto-Legendre quadrature and a time step of $360.0\mathrm{s}$. In the first set of tests the
implicit nonlinear solver was converged to a tolerance of $10^{-14}$, using exact temporal integration of the
variational derivatives \cite{Bauer18,Eldred18,LP21}, allowing for a machine precision error in energy
conservation. In each case the upwinding parameter was set to a value of $\tau=\Delta t/2 = 180.0\mathrm{s}$. 
For the results presented in this section and the next we compute the potential vorticity not through 
the exact time integration of potential enstrophy across the time level as in \eqref{eq::q_linear_in_time} 
or \eqref{eq::q_const_in_time} but rather with a simple second order in time averaging 
of the potential vorticity as $q_h^{n+1/2}=(q_h^n + q_h^{k})/2$ at each Newton iteration $k$, 
where $q_h^n$ and $q_h^k$ are computed instantaneously from \eqref{eq::pv_disc} using
$(\boldsymbol{u}_h^n,h_h^n)$ and $(\boldsymbol{u}_h^k,h_h^k)$ respectively, and without the temporal 
cross terms that would arise from the evaluation of $\overline{q}_h$ via \eqref{eq::q_linear_in_time} 
or \eqref{eq::q_const_in_time}.

\begin{figure}[!hbtp]
\begin{center}
\includegraphics[width=1.0\textwidth,height=0.40\textwidth]{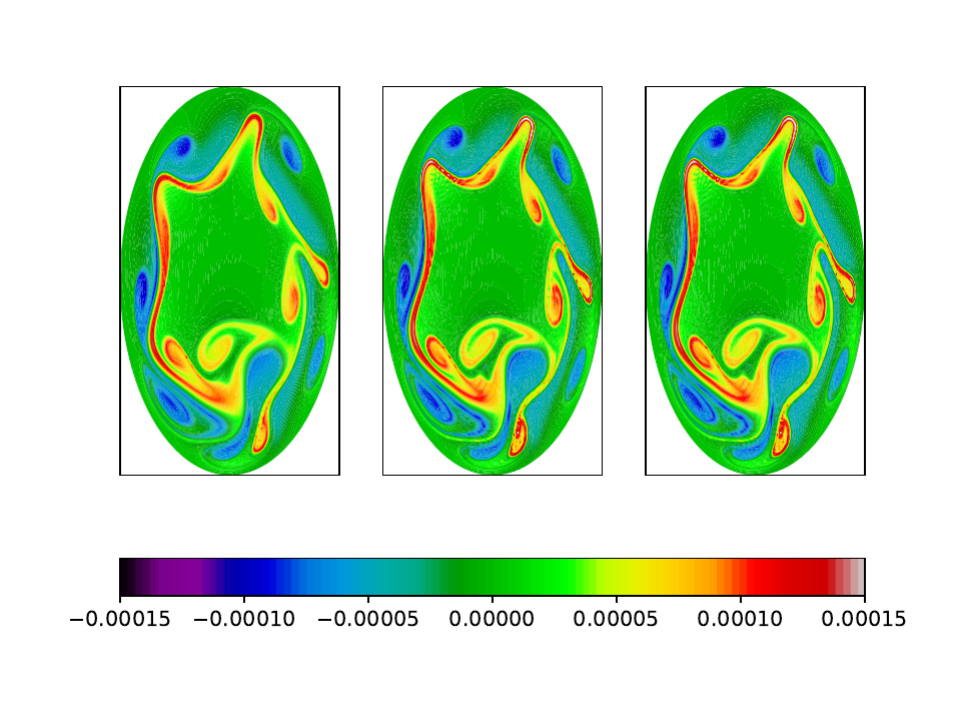}
	\caption{Vorticity field (northern hemisphere) for the Galewsky test case at day 7
	using APVM (left), SUPG (center) and downwinded trial functions (right) for the converged
	solutions.}
\label{fig::vort_nits_day06}
\end{center}
\end{figure}

\begin{figure}[!hbtp]
\begin{center}
\includegraphics[width=1.0\textwidth,height=0.40\textwidth]{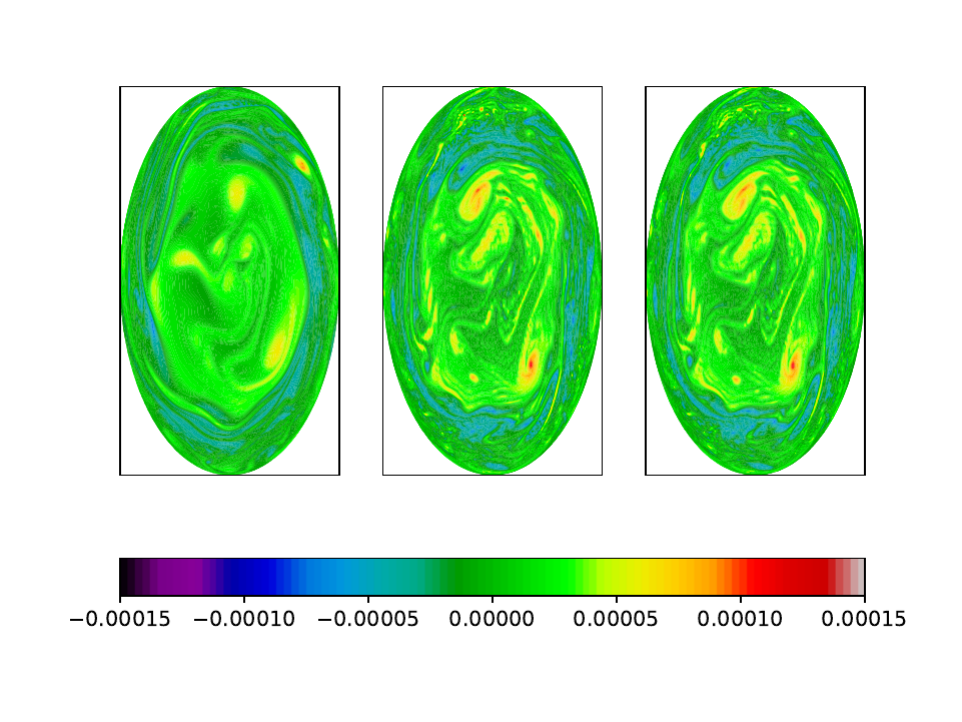}
	\caption{Vorticity field (northern hemisphere) for the Galewsky test case at day 20
	using APVM (left), SUPG (center) and downwinded trial functions (right) for the converged
	solutions.}
\label{fig::vort_nits_day20}
\end{center}
\end{figure}

Figures \ref{fig::vort_nits_day06} and \ref{fig::vort_nits_day20} show the vorticity field in
the northern hemisphere for the three upwinding formulations at days 6 and 20 respectively.
In both cases the SUPG and downwinded trial function schemes look very similar, as suggested by our 
analysis in the previous section. However the APVM
results exhibit excessive dissipation of the potential enstrophy at day 20, due to the symmetric definite
correction to the potential enstrophy evolution that this implies \eqref{eq::pe_adv_apvm}.
The difference between APVM and the other schemes is also evident at day 7. This is due to the
fact that the excessive damping of the potential enstrophy means the nonlinear shear flow instability
is triggered slightly later.

\begin{figure}[!hbtp]
\begin{center}
\includegraphics[width=0.48\textwidth,height=0.32\textwidth]{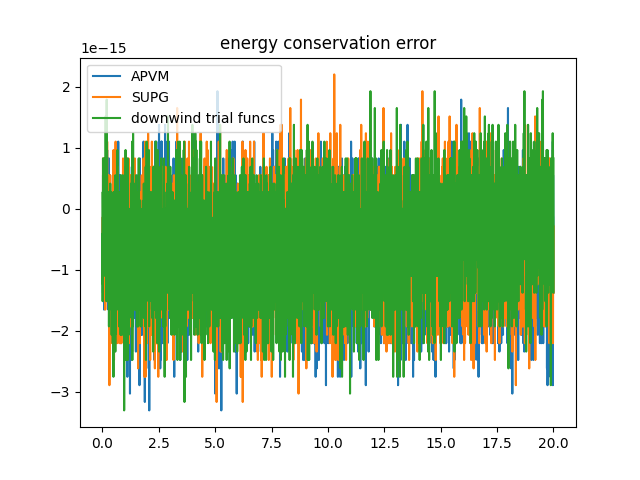}
\includegraphics[width=0.48\textwidth,height=0.32\textwidth]{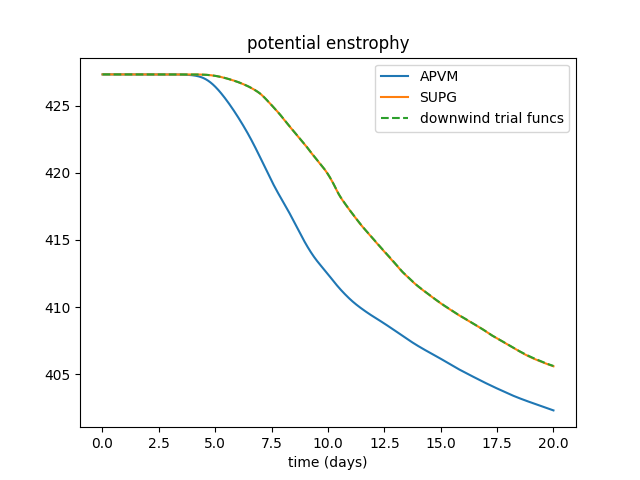}
	\caption{Energy (left) and potential enstrophy (right) conservation for the converged
	solutions.}
\label{fig::conservation_nits}
\end{center}
\end{figure}

The conservation properties (energy and potential enstrophy) are given in Fig. \ref{fig::conservation_nits}.
In all cases the energy conservation error is machine precision. The potential enstrophy
conservation error is comparable for the SUPG and downwinded trial function formulations, reflecting
a similar amount of variation in the vorticity fields, as shown in Fig. \ref{fig::vort_nits_day20}.
However the APVM exhibits a significantly greater amount of potential enstrophy dissipation, owing
to the symmetric definite nature of the correction to the potential enstrophy equation.

\begin{figure}[!hbtp]
\begin{center}
\includegraphics[width=0.48\textwidth,height=0.32\textwidth]{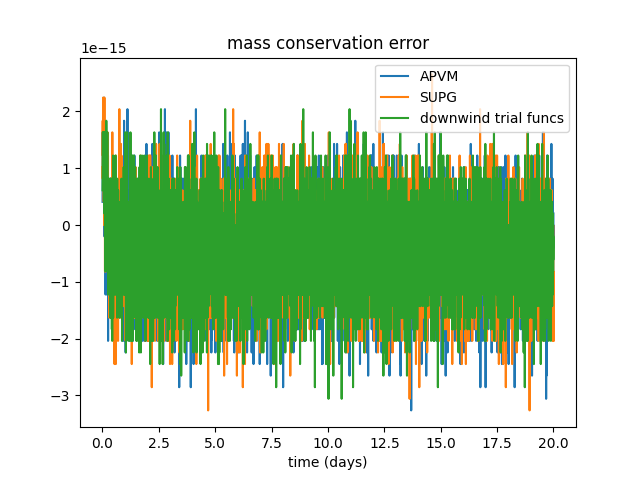}
\includegraphics[width=0.48\textwidth,height=0.32\textwidth]{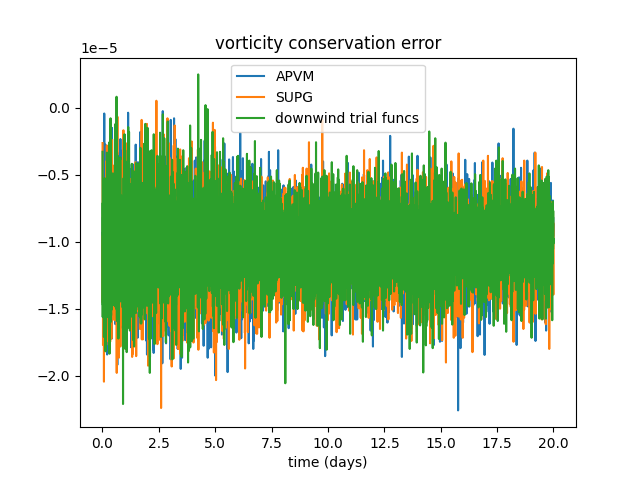}
	\caption{Mass (left) and un-normalised vorticity (right) conservation for the converged
	solutions.}
\label{fig::mass_cons_nits}
\end{center}
\end{figure}

The mass and vorticity conservation errors are given in Fig. \ref{fig::mass_cons_nits}. In each
case the mass conservation is exact (machine precision), while the global vorticity integral
(over a sphere of radius 6371220.0m) is un-normalised and of $\mathcal{O}(10^{-5})$, showing
that in each case the potential vorticity stabilisation has no bearing on its conservation.

\begin{figure}[!hbtp]
\begin{center}
\includegraphics[width=0.48\textwidth,height=0.32\textwidth]{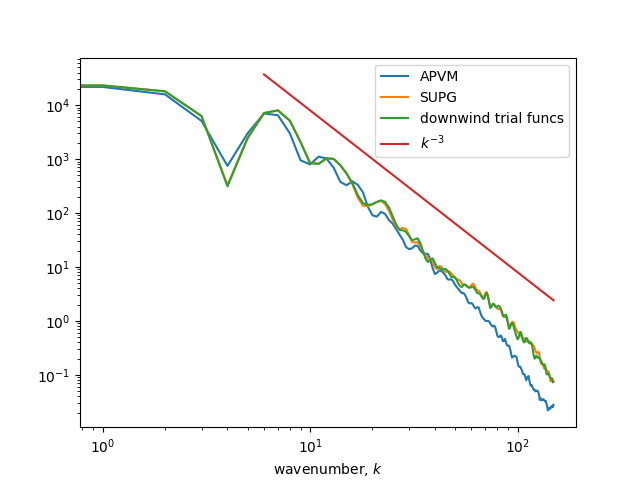}
\includegraphics[width=0.48\textwidth,height=0.32\textwidth]{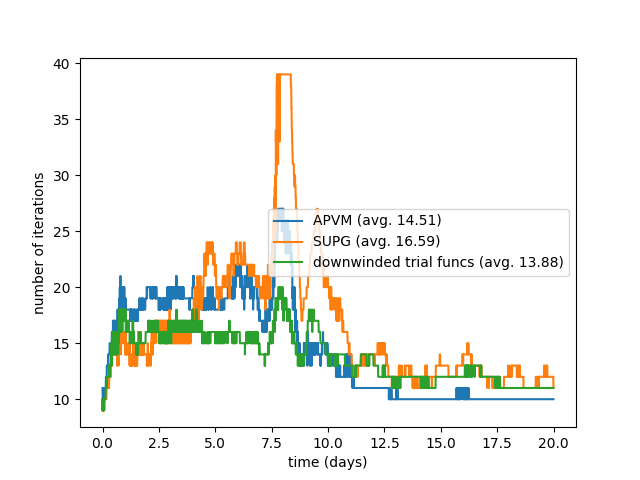}
	\caption{Kinetic energy spectra at day 20 (left), and number of nonlinear iterations
	to achieve convergence (right).}
\label{fig::ke_spectra_nits}
\end{center}
\end{figure}

While the energy and potential enstrophy conservation errors are remarkably similar for the
SUPG and downwinded trial functions, the number of iterations required to achieve convergence
of the nonlinear solver differs greatly, as shown in Fig. \ref{fig::ke_spectra_nits}. The
downwinded trial function solution requires fewer iterations on average to achieve convergence
than either the APVM or SUPG formulations. 
This is perhaps due to the additional higher order corrections in the Taylor series expansion of the
downwinded trial functions \eqref{eq::basis_down}.
This improved solver convergence is despite the fact that the turbulent profile
is almost identical to that for the SUPG formulation, as reflected in the kinetic energy
spectra at day 20, as also shown in Fig. \ref{fig::ke_spectra_nits}, and the potential
vorticity field, as observed in Fig. \ref{fig::vort_nits_day20}. Both the SUPG and downwinded
trial function kinetic energy spectra are closer to the theoretical profile for two dimensional
turbulence of $k^{-3}$ (for spherical harmonic wavenumber $k$).

\begin{figure}[!hbtp]
\begin{center}
\includegraphics[width=0.48\textwidth,height=0.32\textwidth]{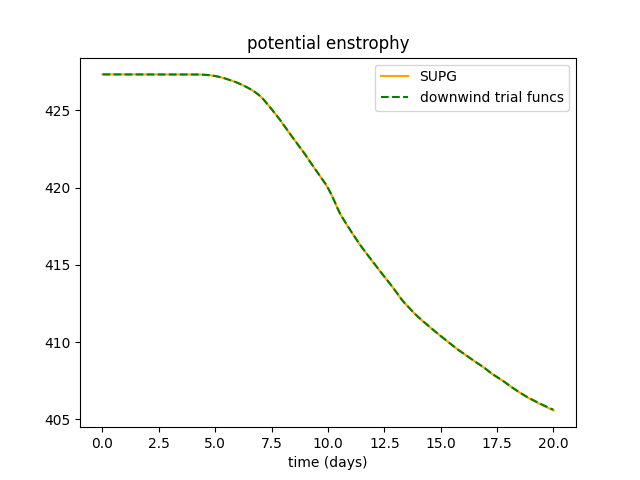}
\includegraphics[width=0.48\textwidth,height=0.32\textwidth]{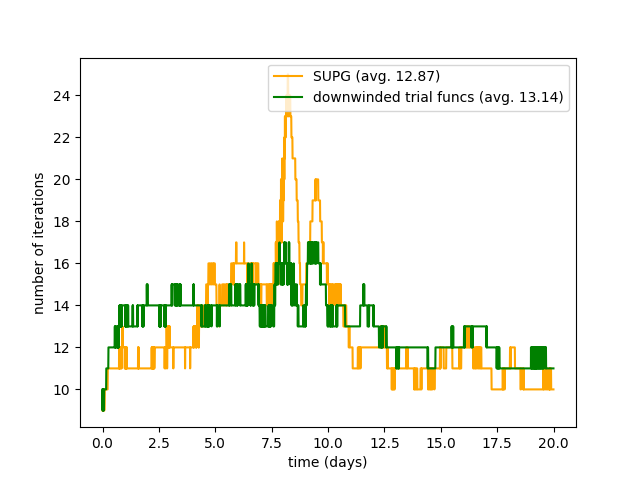}
	\caption{Potential enstrophy conservation (left) and number of iterations
	to achieve convergence (right) using
	$\tau=(2/\Delta t + |\boldsymbol{u}_h|/(2\sqrt{|J|}))^{-1}$.}
\label{fig::tau_dt_udq_nits}
\end{center}
\end{figure}

SUPG schemes are typically run with a more sophisticated representation of $\tau$ so that this
may scale with both the time derivative and the advective term (and any other additional terms
in the underlying equation \cite{AT04}). To this end the experiments were re-run with a value of
$\tau=(2/\Delta t + |\boldsymbol{u}_h|/(2\sqrt{|J|}))^{-1}$. Once more, the potential enstrophy
results were remarkably similar between the SUPG and downwinded trial function schemes. However
the number of iterations required to achieve convergence was reduced for both schemes, with the
SUPG scheme showing the greatest improvement, such that the average number of iterations between
the two schemes was almost the same over the course of the simulations, as show in Fig.
\ref{fig::tau_dt_udq_nits}.

\subsection{Finite number of nonlinear iterations}

The results of the SUPG and downwinded trial function solutions are remarkably similar
for the case where the nonlinear residual is converged to machine precision (at the expense of
a larger number of iterations for the SUPG scheme). However the faster convergence of the downwinded
trial function formulation means that the results differ more sharply for the case where only
a finite number of nonlinear iterations (two) are employed, as is customary in atmospheric models
\cite{Melvin19} for which performance is a premium concern. As for the previous section, here
the potential vorticity is computed as a second order averaging of the instantaneous values at the
previous time level $n$ and the new time level at Newton iteration $k$, $q_h^{n+1/2}=(q_h^n + q_h^{k})/2$, 
and not using the exact temporal integration of the potential enstrophy across the time level, as in
\eqref{eq::q_linear_in_time} or \eqref{eq::q_const_in_time}.

\begin{figure}[!hbtp]
\begin{center}
\includegraphics[width=1.0\textwidth,height=0.40\textwidth]{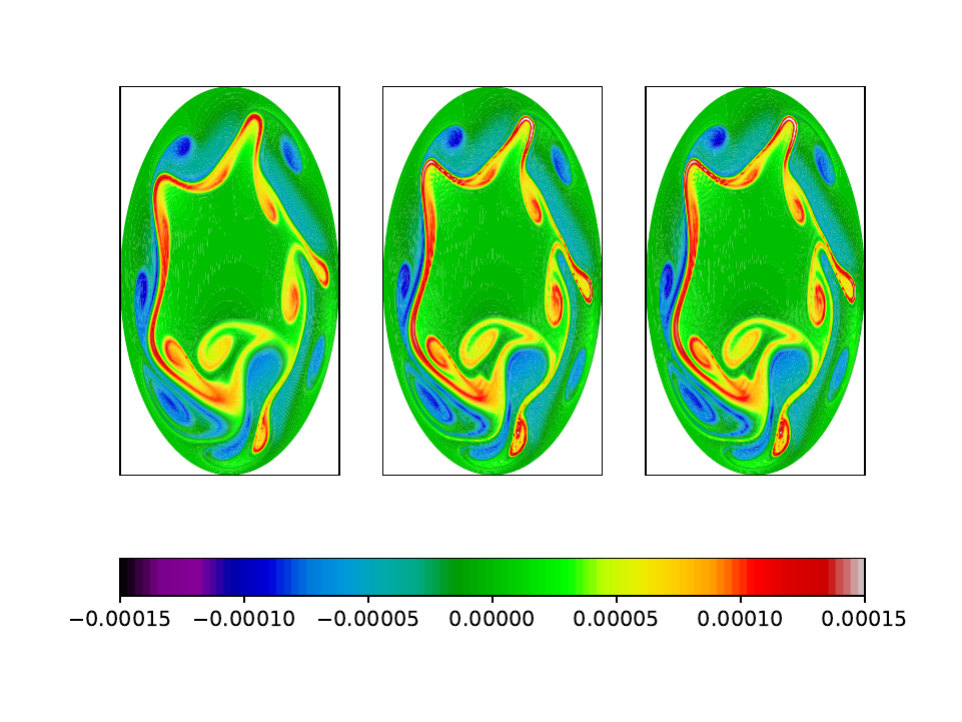}
	\caption{Vorticity field (northern hemisphere) for the Galewsky test case at day 7
	using APVM (left), SUPG (center) and downwinded trial functions (right) using two
	nonlinear iterations per time step.}
\label{fig::vort_2its_day06}
\end{center}
\end{figure}

\begin{figure}[!hbtp]
\begin{center}
\includegraphics[width=1.0\textwidth,height=0.40\textwidth]{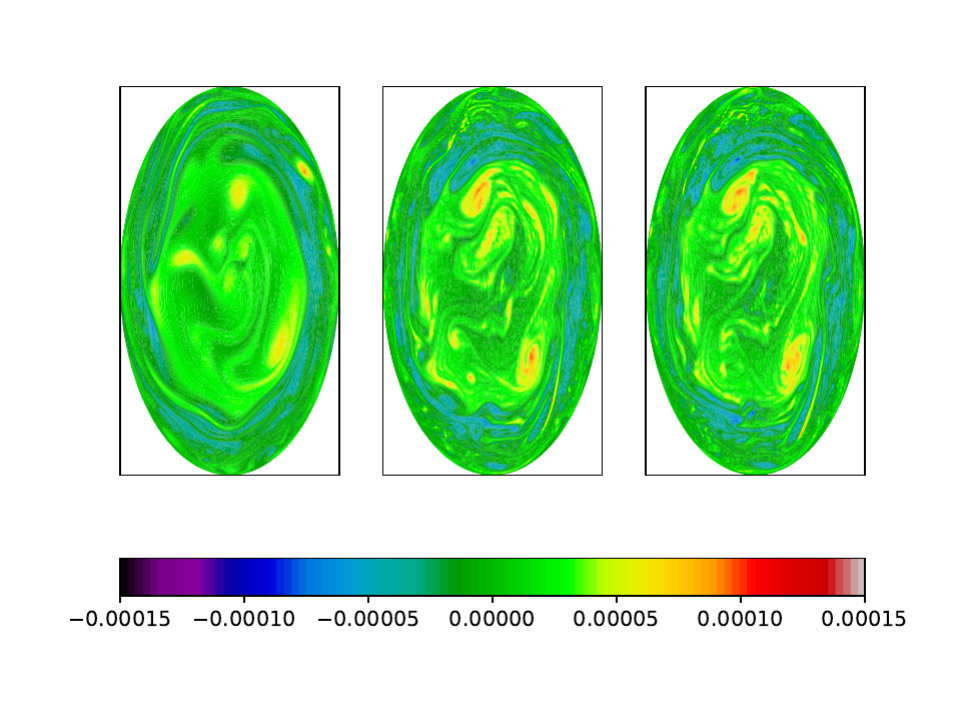}
	\caption{Vorticity field (northern hemisphere) for the Galewsky test case at day 20
	using APVM (left), SUPG (center) and downwinded trial functions (right) using two
	nonlinear iterations per time step.}
\label{fig::vort_2its_day20}
\end{center}
\end{figure}

\begin{figure}[!hbtp]
\begin{center}
\includegraphics[width=0.48\textwidth,height=0.32\textwidth]{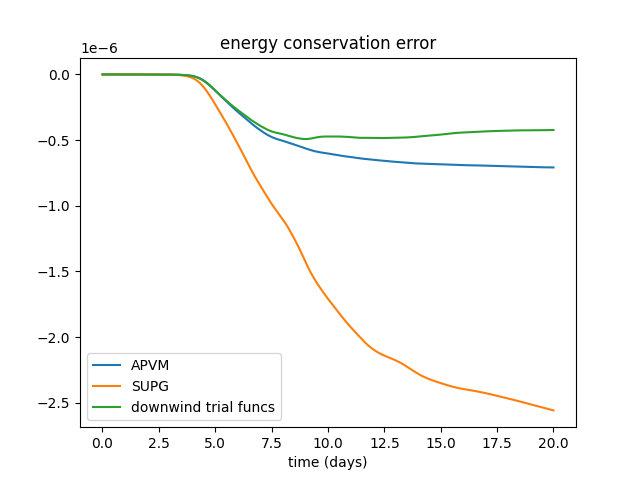}
\includegraphics[width=0.48\textwidth,height=0.32\textwidth]{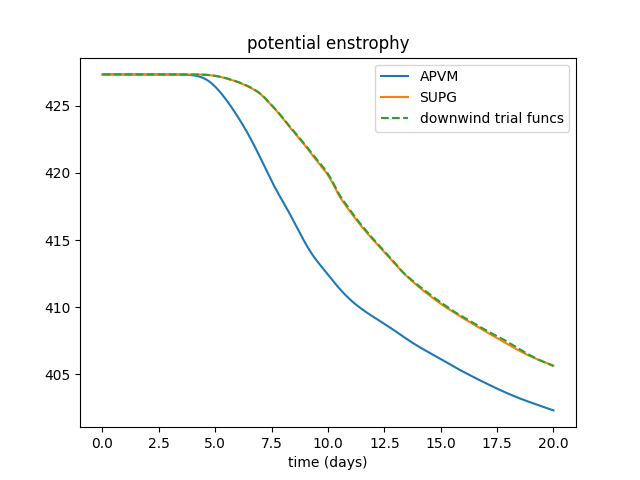}
	\caption{Energy (left) and potential enstrophy (right) conservation using two
	nonlinear iterations per time step.}
\label{fig::conservation_2its}
\end{center}
\end{figure}

\begin{figure}[!hbtp]
\begin{center}
\includegraphics[width=0.48\textwidth,height=0.32\textwidth]{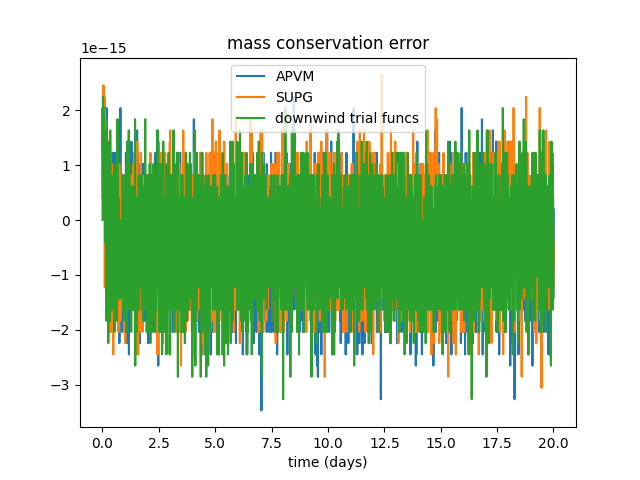}
\includegraphics[width=0.48\textwidth,height=0.32\textwidth]{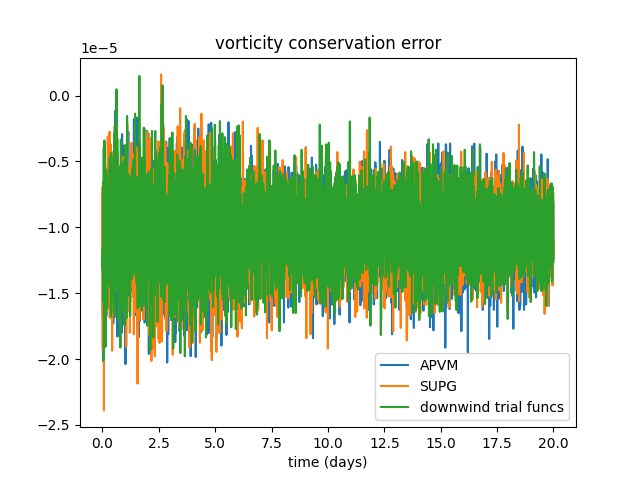}
	\caption{Mass (left) and un-normalised vorticity (right) conservation using two
	nonlinear iterations per time step.}
\label{fig::mass_cons_2its}
\end{center}
\end{figure}

As shown in Figs. \ref{fig::vort_2its_day06} and \ref{fig::vort_2its_day20}, the results of the
SUPG and downwinded trial function formulations are still remarkably similar to visual inspection
at both days 6 and 20 for the case where only two nonlinear iterations are employed per time
step. While the potential enstrophy errors are also very similar, as observed in Fig.
\ref{fig::conservation_2its}, the energy conservation errors differ markedly, with the downwinded
trial function energy conservation error being significantly smaller over the course of the
simulation. The mass and vorticity conservation errors are essentially the same as for the
converged solution, as shown in Fig. \ref{fig::mass_cons_2its}.

\begin{figure}[!hbtp]
\begin{center}
\includegraphics[width=0.48\textwidth,height=0.32\textwidth]{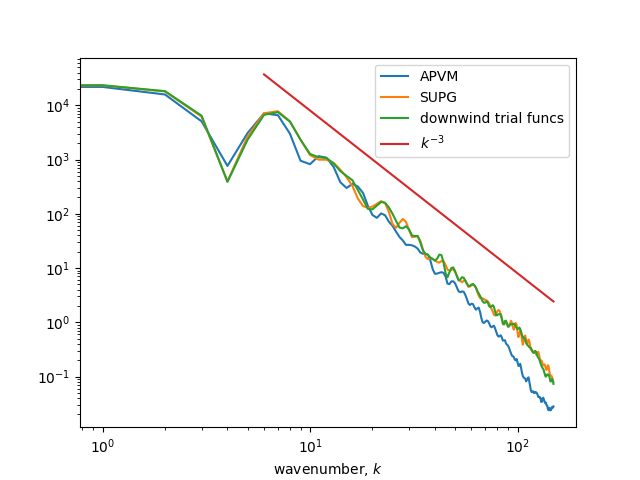}
\includegraphics[width=0.48\textwidth,height=0.32\textwidth]{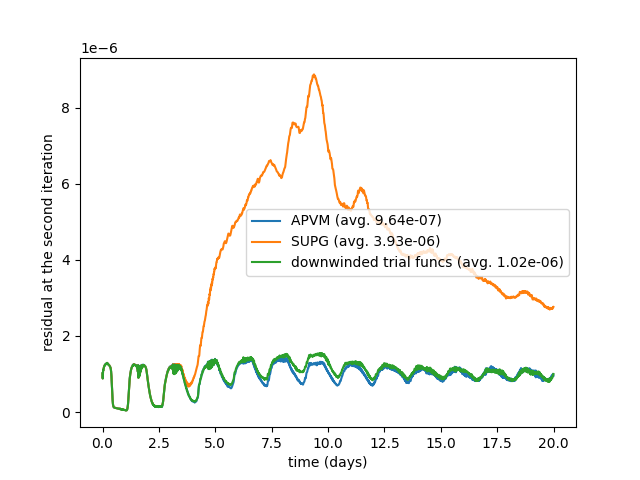}
	\caption{Kinetic energy spectra at day 20 (left), and residual error at the second
	iteration using two nonlinear iterations per time step (left).}
\label{fig::ke_spectra_2its}
\end{center}
\end{figure}

This difference in the energy conservation error is also reflected in the nonlinear
residual at the second iteration, as shown in Fig. \ref{fig::ke_spectra_2its}. This error
is significantly smaller for the downwinded trial functions than the SUPG formulation, particularly
at the maximum growth of the shear instability (around days 8-9). Curiously, the residual error
is remarkably similar for the APVM and downwinded trial functions, despite the fact that the
turbulence profile is much richer for the downwinded trial functions, as also shown for day 20
in Fig. \ref{fig::ke_spectra_2its}.

\begin{figure}[!hbtp]
\begin{center}
\includegraphics[width=0.48\textwidth,height=0.32\textwidth]{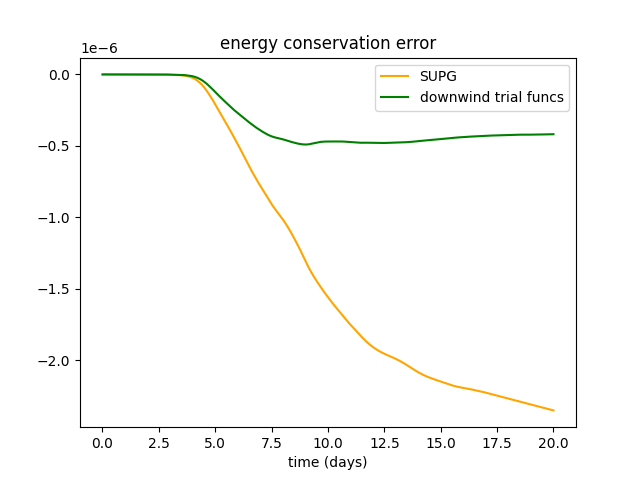}
\includegraphics[width=0.48\textwidth,height=0.32\textwidth]{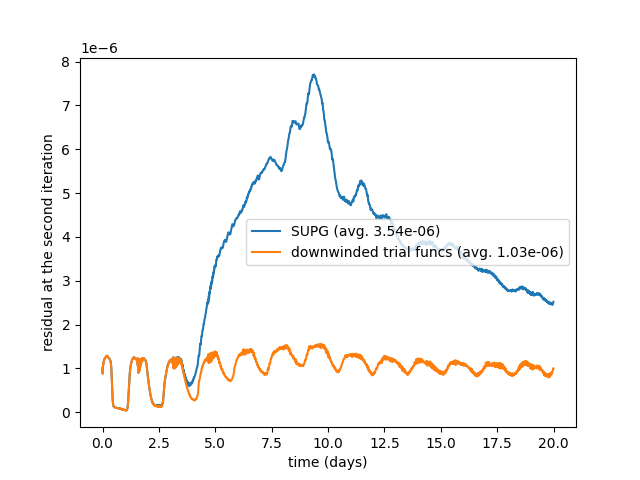}
	\caption{Energy conservation (left) and residual errors
		 at the second iteration (right), using two nonlinear iterations
	and $\tau=(2/\Delta t + |\boldsymbol{u}_h|/(2\sqrt{|J|}))^{-1}$.}
\label{fig::tau_dt_udq_2its}
\end{center}
\end{figure}

The energy conservation and residual errors between the SUPG and downwinded trial function formulations
are also compared using a spatially and temporally varying value of
$\tau=(2/\Delta t + |\boldsymbol{u}_h|/(2\sqrt{|J|}))^{-1}$ in Fig. \ref{fig::tau_dt_udq_2its}.
While the residual errors for the SUPG formulation are somewhat improved with respect to the
constant value results given in Fig. \ref{fig::ke_spectra_2its}, the energy energy conservation
errors are significantly poorer than those shown for the SUPG scheme in Fig. \ref{fig::conservation_2its}.
The downwinded trial function results are similar in both cases.

\subsection{Exact temporal integration of potential enstrophy}

\begin{figure}[!hbtp]
\begin{center}
\includegraphics[width=0.48\textwidth,height=0.32\textwidth]{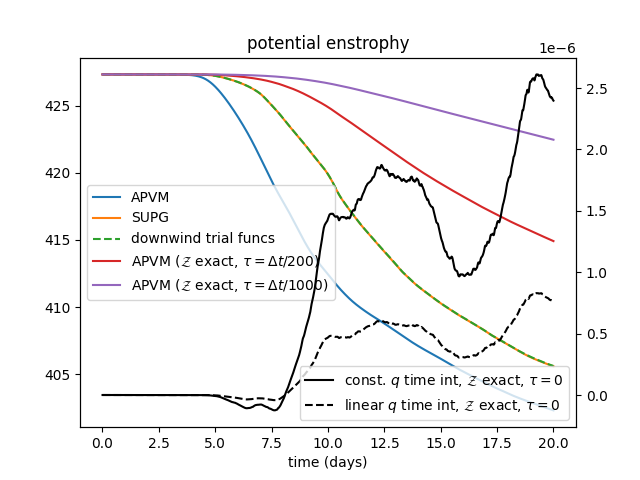}
\includegraphics[width=0.48\textwidth,height=0.32\textwidth]{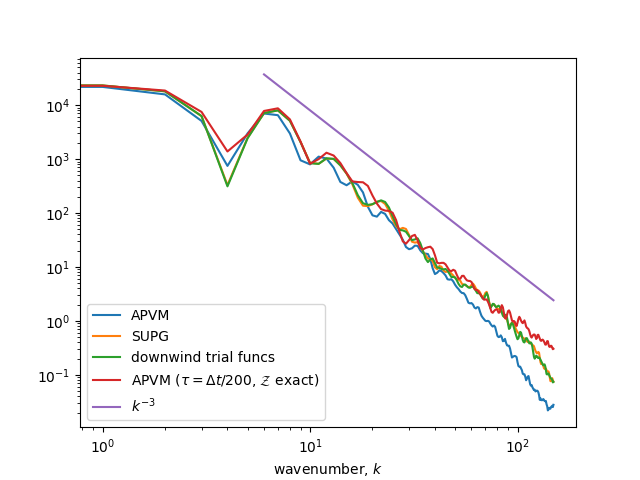}
	\caption{Potential enstrophy (left) and kinetic energy spectra at day 20 (right)
using both inexact and exact potential enstrophy conserving time stepping.
	Note the different scales of the vertical axes for the stabilised configuration absolute 
	values vs. unstabilised configuration normalised errors in the left hand plot.}
\label{fig::pot_enst_con}
\end{center}
\end{figure}

We also compare the energy conserving results presented in Section 4.1 with results for which potential
enstrophy is exactly integrated across each time level, via the diagnosis of $\overline{q}_h$ as
either piecewise linear \eqref{eq::q_linear_in_time} or piecewise constant \eqref{eq::q_const_in_time} in time. 
Since there is no internal dissipation of energy in space or time, and only minimal change in potential enstrophy 
for these configurations, any injection of potential
enstrophy, such as that provided by the SUPG or downwinded trial function stabilisation schemes, will lead to
instability. Consequently these formulation can only be run in either a neutral state, with no damping, or in
conjunction with a symmetric positive definite damping term, such as that for the APVM. Results showing the potential
enstrophy evolution for the exact potential enstrophy integration formulations are presented in Fig.
\ref{fig::pot_enst_con}, with these contrasted against those for the standard second order representation as
shown in Fig. \ref{fig::conservation_nits}. As observed, the second order in time representation of potential
vorticity lead to a potential enstrophy conservation error approximate one quarter that of the piecewise 
constant representation for the completely undamped case ($\tau=0$). As previously noted, since these schemes
ensure the exact temporal integration of $\mathcal{Z}_h$ across the time level, but not between time levels,
while they yield greatly improved results in terms of stability and conservation, unlike energy $\mathcal{Z}_h$
is still not conserved exactly.

As observed in Fig. \ref{fig::pot_enst_con}, using the exact temporal integration of the potential enstrophy
the APVM may be run stably with any choice of $\tau\ge 0$, including smaller values than those required to 
stabilise simulations with an inexact time integration of potential enstrophy. Using such small values
of $\tau$, the APVM may yield superior results to either the SUPG or downwinded trial function formulations in terms
of both conservation error and turbulent spectra. 

\subsection{Comparisons in a gravity wave dominated regime}

While the previous test provides useful insights into the performance of the different potential vorticity
upwinding schemes for flows dominated by vorticity dynamics and turbulence, it is less useful for investigating
the different schemes in the context of gravity wave dominated flows. In order to compare the different schemes
in this regime we apply them to a standard test case for orographically generated gravity waves \cite{Williamson92}.
The test involves a geostrophically balanced flow for a fluid of mean depth $H=5960m$ over an isolated 
mountain with a profile of $b(r)=2000\times(1-r/R)$, where $R=\pi/9$ and $r^2=(\theta-\pi/6)^2 + (\lambda+\pi/2)^2$,
and $\theta$ and $\lambda$ are the latitude and longitude respectively. The bottom topography is applied to 
the solver as an addition to the Bernoulli potential, $\overline{P}_h$ \eqref{eq::dH_dh} of the form $\int g\phi_hb_h\mathrm{d}\Omega$,
where $b_h$ is the discrete projection of $b$ onto $\mathbb{V}_2$. The topography is correspondingly incorporated
into the energy as
$\mathcal{H}_h = \int\frac{1}{2}h_h\boldsymbol{u}_h\cdot\boldsymbol{u}_h + \frac{1}{2}g(h_h+b_h)^2\mathrm{d}\Omega$
\cite{Wimmer20}.

\begin{figure}[!hbtp]
\begin{center}
\includegraphics[width=0.48\textwidth,height=0.32\textwidth]{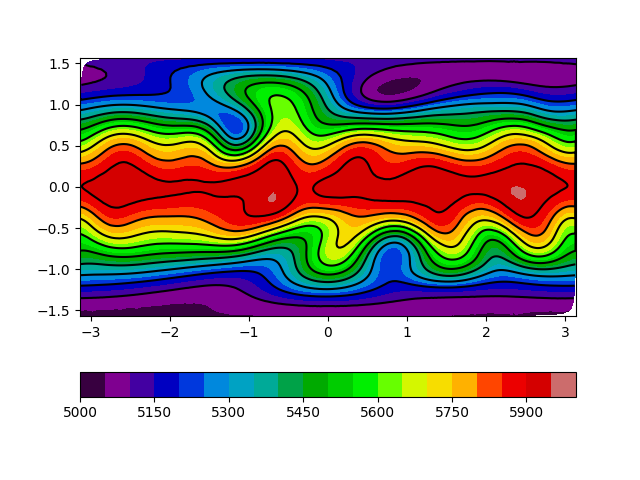}
\includegraphics[width=0.48\textwidth,height=0.32\textwidth]{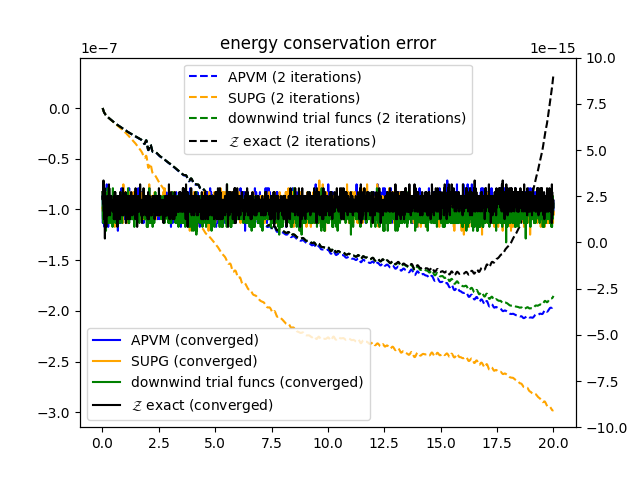}
	\caption{Left: total depth at day 15 for the nonlinearly converged solution of the 
	flow over an isolated mountain test case at each 
	time step using the downwinded trial function stabilisation. Solid contours
	are given from $5100m$ to $5900m$ in increments of $100m$.
	Right: normalised total energy
	conservation errors using two Newton iterations and nonlinear converged
	solutions for the different potential vorticity upwinding schemes. Note
	the different vertical axes for the energy conservation error using two
	iterations and nonlinearly converged solutions.}
\label{fig::williamson5_1}
\end{center}
\end{figure}

The different schemes were all run using $6\times 16\times 16$ elements using polynomials of degree $p=3$ and 
8 point quadrature with a time step of $\Delta t=600s$ for a total of 20 days. As for the barotropic jet test 
configuration above, these were compared for
simulations that were run to nonlinear convergence at each time step, and also for simulations using just two 
Newton iterations per time step. Visual inspection of the total depth ($h_h + b_h$) and vorticity fields at day 15 showed little 
variation between the result for the different formulations, with the nonlinearly converged solutions being in
all cases slightly sharper than the results using just two nonlinear iterations. The total depth
field at day 15 using the downwinded trial function stabilisation for the nonlinearly converged 
solution is presented in Fig. \ref{fig::williamson5_1}.

While the results appear visually similar, differences in the energy
conservation errors can be easily discerned in Fig. \ref{fig::williamson5_1}. As expected, in
all cases energy is conserved to machine precision when the Newton solver is run to convergence.
When using only two Newton iterations per time step the SUPG method leads to the largest energy
loss, while the exact computation of the potential enstrophy across the time level, via 
\eqref{eq::q_const_in_time}, leads to long term energy growth. As for the previous test case all formulations 
other than the exact integral of potential enstrophy across the time level use a simple time centered formulation
of $q_h^{n+1/2} = (q_h^n + q_h^k)/2$, with $q_h^n$ and $q_h^k$ computed instantaneously via \eqref{eq::pv_disc}.
Long time energy growth is also slightly perceptible for the APVM and downwinded trial function 
solutions, suggesting that these methods are perhaps slightly under-damped.

\begin{figure}[!hbtp]
\begin{center}
\includegraphics[width=0.48\textwidth,height=0.32\textwidth]{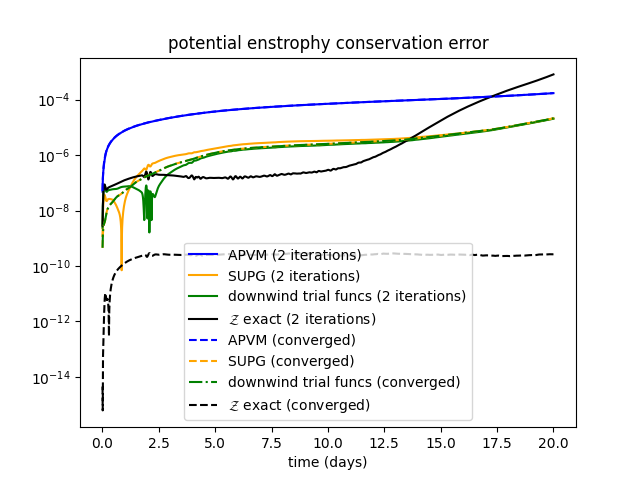}
\includegraphics[width=0.48\textwidth,height=0.32\textwidth]{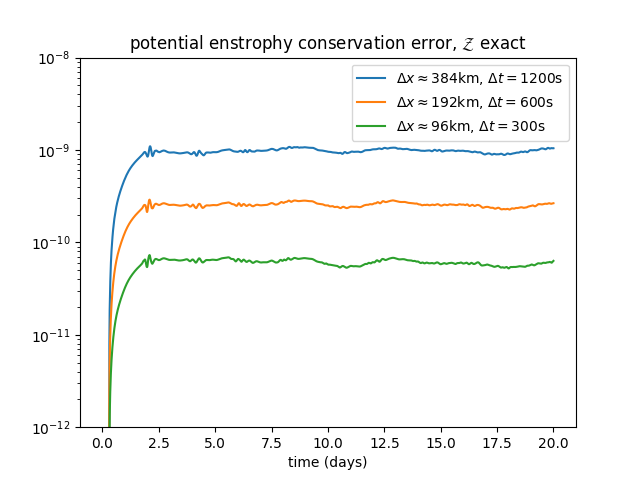}
	\caption{Left: normalised absolute values of the potential enstrophy
	conservation errors using two Newton iterations and nonlinear converged
	solutions for the different potential vorticity upwinding schemes for the flow over 
	an isolated mountain test case. Note that
	these are decaying for all schemes, with the exception of the two iteration 
	$\mathcal{Z}$ exact solution. Right: normalised absolute value of the potential enstrophy
	conservation error for the nonlinear converged
	solution as a function of resolution for the piecewise constant in time 
	exact integration of $\mathcal{Z}$ across the time level via \eqref{eq::q_const_in_time}.}
\label{fig::williamson5_2}
\end{center}
\end{figure}

The potential enstrophy conservation errors are presented in Fig. 
\ref{fig::williamson5_2}, in normalised absolute value form on a logarithmic plot. All solutions
exhibit a decay of potential enstrophy with time, with the exception of the two iteration solution
using exact integration of potential enstrophy across the time level, for which this grows with
time, suggesting that for two iterations the solution will ultimately go unstable. Once more the
APVM method is most dissipative, using both two Newton iterations or a nonlinearly converged
solution. Results are nearly identical for the SUPG and downwinded trial function solutions for
nonlinearly converged solutions, while for two iterations the SUPG method is slightly more
dissipative. As for the previous test case, exactly integrating the potential enstrophy across 
the time level with a converged Newton iteration at each time step leads to
greatly improved potential enstrophy conservation errors, which are steady throughout the
simulation at $\mathcal{O}(10^{-10})$.

\begin{figure}[!hbtp]
\begin{center}
\includegraphics[width=0.48\textwidth,height=0.32\textwidth]{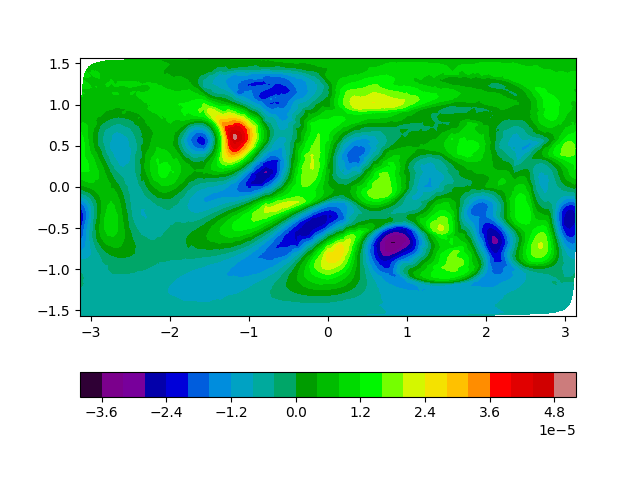}
\includegraphics[width=0.48\textwidth,height=0.32\textwidth]{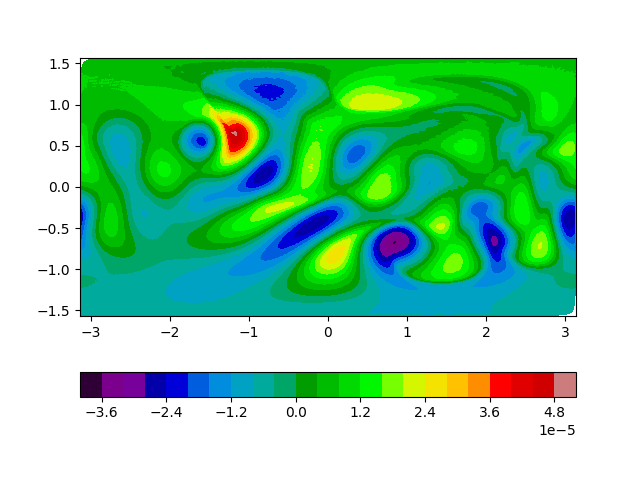}
	\caption{
		Vorticity field at day 15 using $6\times 8^2$ elements and $\Delta t=1200$s (left) and
		$6\times 32^2$ elements and $\Delta t=300$s (right) for the flow over an isolated 
		mountain test case using piecewise constant in time exact integration of $\mathcal{Z}$ 
		across the time level.}
\label{fig::williamson5_3}
\end{center}
\end{figure}

Since the potential enstrophy shows no discernible decrease with time when this is exactly integrated across
the time level, its turbulent cascade to grid scales will introduce some amount of noise into the 
solution. In Fig. \ref{fig::williamson5_2} we give the normalised potential enstrophy conservation 
errors as a function of time with increased spatial and temporal resolution. These errors decrease 
quadratically (the order of the time stepping scheme), suggesting that while there will be some
degree of noise for these solutions due to the turbulent cascade of potential enstrophy in the absence
of dissipation, this noise nevertheless decreases with resolution. This result is visually confirmed 
in Fig. \ref{fig::williamson5_3}, which gives the vorticity field at day 15 for the flow
over the isolated mountain test case using resolutions of $6\times 8^2$ elements with $\Delta t=1200s$
and $6\times32^2$ elements with $\Delta t=300s$, where it is clearly observed that the higher 
resolution result is also smoother at the grid scale. As an aside, we also note that if the 
potential enstrophy does not decrease with time for two dimensional flows where this cascades to the grid scale,
this is most likely indicative of grid scale noise, which may be deemed unacceptable for practical applications.

As for the previous test case, mass and vorticity are exactly conserved in all cases. These 
results are not presented here. In terms of the total number of Newton iterations for the 
case where the different configurations were run to convergence at each time step, these
were remarkably similar, averaging 23.27 iterations for the APVM scheme, and 23.26 iterations
for the other formulations. The results were also very similar between the different formulations
for the residual error at the second iteration when run with just two Newton iterations. The
close similarity of these results is perhaps unsurprising, since the upwind stabilisation 
applied here is specifically targeting the potential vorticity transport, and this has less 
overall effect on the dynamics for the gravity wave dominated test case than it does for the 
previous shear instability test case.

\section{Conclusions}

Different forms of potential vorticity upwinding are analysed and compared for a high order mixed finite element 
discretisation of the rotating shallow water equations on the sphere. These include the well known APVM and SUPG 
methods, as well as a more recently proposed method by which the basis functions are evaluated at downstream 
locations in reference element coordinates. Both analysis and numerical experiments show that the reference 
element downwinding and SUPG methods are clearly superior to the APVM method in terms of their turbulent 
profiles and potential enstrophy conservation. While the SUPG and downwinded trial function formulations are 
very similar in terms of both their potential enstrophy conservation and turbulent spectra, the downwinded trial 
functions exhibit better energy conservation and smaller residual errors for a finite number of nonlinear 
iterations for both constant and non-constant values of the upwinding parameter. For a fully converged solution 
the downwinded trial functions exhibit faster convergence if a constant upwinding parameter is used, however the 
convergence rates are more comparable for a non-constant upwinding parameter.

One possible explanation for the improved results of the downwinded trial functions is that while the potential
vorticity is diagnosed from the prognostic variables ($\boldsymbol{u}_h, h_h$) in the formulation presented here,
it may alternatively be represented via a material advection equation \eqref{eq::pv_adv}. By applying a higher
order representation of the upwinding, we are effectively providing a more accurate approximation of where the
potential vorticity departure point would be located at some previous time, were this to be represented as a
Lagrangian (or semi-Lagrangian) variable. The downwinded trial functions therefore take more account of the
material nature of the potential vorticity.

In addition to the comparison of upwind stabilisation schemes, we also present new temporal formulations of the
potential vorticity diagnostic equation that imply the exact integration of potential enstrophy across the time 
level (but not between time levels) for both piecewise linear and piecewise constant (in time) potential vorticity.
These temporal discretisations are observed to be almost neutrally stable, and so allow for the stable simulation 
of a mature turbulent state without any form of dissipation. When run in conjunction with the APVM using values of 
$\tau$ $\mathcal{O}(100)$ times smaller than those typically required, these formulations exhibit superior 
conservation properties and turbulence profiles than those observed for either the SUPG or downwinded trial 
function formulations when run using a simple time centered (inexact) formulation for the potential vorticity and 
a more moderate value of the upwinding parameter.

Results for a standard test case for flow over an isolated mountain (for which the solution is well resolved 
at lower resolution than for the geostrophic turbulence test case) using exact integration of potential enstrophy
across the time level show that potential enstrophy conservation errors decrease with the second order accuracy of 
the time stepping scheme. As such the grid scale noise due to the lack of potential enstrophy damping in this 
formulation reduces with increased resolution, which is also confirmed by visual inspection of the vorticity field
at different resolutions.

We end by noting that while the upwind stabilisation of the potential vorticity has been the main subject of this
article, these approaches are equally applicable to other material variables that appear in the skew-symmetric
operators of Hamiltonian systems, as has been previously demonstrated for the potential temperature in the case
of the 3D compressible Euler equations \cite{LP21,WCB21}.

\section*{Acknowledgments}

\thethanks



\begin{thebibliography}{10}
\expandafter\ifx\csname url\endcsname\relax
  \def\url#1{\texttt{#1}}\fi
\expandafter\ifx\csname urlprefix\endcsname\relax\def\urlprefix{URL }\fi
\expandafter\ifx\csname href\endcsname\relax
  \def\href#1#2{#2} \def\path#1{#1}\fi

\bibitem{SB85}
R.~Sadourny, C.~Basdevant, Parameterization of {S}ubgrid {S}cale {B}arotropic and {B}aroclinic {E}ddies
  in {Q}uasi-geostrophic {M}odels: {A}nticipated {P}otential {V}orticity {M}ethod, J. Atmos. Sci. 42
  (1985) 1353--1363

\bibitem{BH82}
A.~N. Brooks, T.~J.~R.~Hughes, Streamline upwind/Petrov-Galerkin formulations for convection dominated flows
with particular emphasis on the incompressible Navier-Stokes equations, Comput. Meth. Appl. Mech. Engrg. 32
(1982) 199--259

\bibitem{Hughes98}
T.~J.~R.~Hughes, G.~R.~Feijoo, L.~Mazzei, J.-B.~Quincy, The variational multiscale method--a paradigm for
computational mechanics, Comput. Methods Appl. Mech. Engrg. 166 (1998) 3--24

\bibitem{MC14}
A.~T.~T.~McRae, C.~J.~Cotter, Energy- and enstrophy conserving schemes for the
  shallow water equations, based on mimetic finite elements, Q. J. R. Meteorol.
  Soc. 140 (2014) 2223--2234.

\bibitem{NC17}
A.~Natale, C.~J. Cotter, Scale-selective dissipation in energy-conserving finite-element schemes for 
  two-dimensional turbulence, Q. J. R. Meteorol. Soc. 143 (2017) 1734--1745.

\bibitem{Natale16}
A.~Natale, J.~Shipton, C.~J. Cotter, Compatible finite element spaces for
  geophysical fluid dynamics, Dyn. Stat. Climate Sys. 1 (2016) 1--31.

\bibitem{Bauer18}
W.~Bauer, C.~J. Cotter, Energy-enstrophy conserving compatible finite element
  schemes for the rotating shallow water equations with slip boundary
  conditions, J. Comp. Phys. 373 (2018) 171--187.

\bibitem{WCB21}
G.~A.~Wimmer, C.J.~Cotter, W.~Bauer, Energy conserving {SUPG} methods for
  compatible finite element schemes in numerical weather prediction,
  SMAI journal of computational mathematics, 7 (2021), 267--300

\bibitem{Brecht21}
R.~Brecht, W.~Bauer, A.~Bihlo, F.~Gay-Balmaz, S.~MacLachlan, Selective decay for the rotating shallow-water 
  equations with a structure-preserving discretization, Physics of Fluids 33 (2021) 116604.

\bibitem{MG15}
S.~Marras, F.~X.~Giraldo, A parameter-free dynamic alternative to hyper-viscosity
  for coupled transport equations: {A}pplication to the simulation of 3{D} squall
  lines using spectral elements, J. Comp. Phys. 283 (2015) 360--373.

\bibitem{Lee21}
D.~Lee, Petrov-Galerkin flux upwinding for mixed mimetic spectral elements, and
  its application to geophysical flow problems, Comput. Math. Appl. 89 (2021) 68--77

\bibitem{LP21}
D.~Lee,  A.~Palha, Exact spatial and temporal balance of energy exchanges within a
  horizontally explicit/vertically implicit non-hydrostatic atmosphere, J. Comp. Phys.
  440 (2021) 110432

\bibitem{Galewsky04}
J.~Galewsky, R.~K.~Scott, L.~M.~Polvani, An initial-value problem for testing numerical models of the
  global shallow water equations, Tellus 56A (2004) 429--440

\bibitem{LPG18}
D.~Lee, A.~Palha, M.~Gerritsma, Discrete conservation properties for shallow
  water flows using mixed mimetic spectral elements, J. Comp. Phys. 357 (2018)
  282--304

\bibitem{Eldred18}
C.~Eldred, T.~Dubos, E.~Kritsikis, A quasi-{H}amiltonian discretization of the
  thermal shallow water equations, J. Comp. Phys. 379 (2019) 1--31

\bibitem{Celledoni12}
E.~Celledoni, V.~Grimm, R.~I.~McLachlan, D.~I.~McLaren, D.~O'Neale, B.~Owren, 
  G.~R.~W.~Quispel, Preserving energy resp. dissipation in numerical {PDE}s using
  the "average vector ﬁeld" method, J. Comput. Phys. 231 (2012) 6770--6789.

\bibitem{Hairer11}
D.~Cohen, E.~Hairer, Linear energy-preserving integrators for Poisson systems,
  BIT Numer. Math. 51 (1) (2011) 91–101.

\bibitem{Lee21b}
D.~Lee, An energetically balanced, quasi-{N}ewton integrator for non-hydrostatic
  vertical atmospheric dynamics, J. Comp. Phys. (2021) 109988.

\bibitem{Wimmer20}
G.~A.~Wimmer, C.J.~Cotter, W.~Bauer, Energy conserving upwinded compatible
  finite element schemes for the rotating shallow water equations, J. Comp. Phys.
  401 (2020) 109016

\bibitem{AT04}
J.~Ed Akin, T.~E.~Tezduyar, Calculation of the advective limit of the {SUPG}
  stabilization parameter for linear and higher-order elements, Comput. Methods
  Appl. Mech. Engrg. 193 (2004) 1909--1922

\bibitem{Rognes06}
M.~E.~Rognes, D.~A.~Ham, C.~J.~Cotter, A.~T.~T.~McRae, Automating the solution of
  {PDE}s on the sphere and other manifolds in {FE}ni{CS} 1.2, Geosci. Model. Dev.
  6 (2013) 2099--2119

\bibitem{LP18}
D.~Lee, A.~Palha, A mixed mimetic spectral element model of the rotating
  shallow water equations on the cubed sphere, J. Comp. Phys. 375 (2018)
  240--262

\bibitem{Melvin19}
T.~Melvin, T.~Benacchio, B.~Shipway, N.~Wood, J.~Thuburn, C.~Cotter, A mixed
  finite-element, finite-volume, semi-implicit discretisation for atmospheric
  dynamics: {C}artesian geometry, Q. J. R. Meteorol. Soc. (2019) 1--19.

\bibitem{Williamson92}
D.~L.~Williamson, J.~B.~Drake, J.~J.~Hack, R.~Jakob, P.~N.~Swarztrauber, 
  A standard test set for numerical approximations to the shallow water 
  equations in spherical geometry, J. Comp. Phys. 102, (1992) 211--224.

\end{thebibliography}
\end{document}